# JACOBI FORMS AND THE STRUCTURE OF DONALDSON INVARIANTS FOR 4-MANIFOLDS WITH $b_+ = 1$

LOTHAR GÖTTSCHE AND DON ZAGIER

## Contents

1. INTRODUCTION

In the last few years the Donaldson invariants of smooth 4-manifolds have been very powerful tools. They are defined using a Riemannian metric $g$ on the 4-manifold $X$, but in the case that $b_+(X) > 1$ they are independent of $g$. Recently very much progress has been made in the understanding of these invariants. In particular Kronheimer and Mrowka [Kr-M1], [Kr-M2] and Fintushel and Stern [F-S2] have shown some important structure theorems.

For $C \in H^2(X, \mathbb{Z})$ and $c_2 \in \mathbb{Z}$ we denote by $\Phi_{C,d}^X$ the Donaldson invariant of $X$ with respect to a complex rank 2 bundle $\mathcal{E}$ with Chern classes $C$ and $c_2$ with $4c_2 - C^2 - \frac{3}{2}(1 + b_+(X)) = d$. $\Phi_{C,d}^X$ is a linear map from the set $A_d(X)$ of polynomials of degree $d$ in $H_2(X, \mathbb{Q}) \oplus H_0(X, \mathbb{Q})$ to $\mathbb{Q}$, where the classes in $H_2(X, \mathbb{Q})$ are given degree 1 and the class $p \in H_0(X, \mathbb{Z})$ corresponding to a point is given degree 2. We put $\Phi_{C,d}^X := 0$ if $b_+(X)$ is even or if $d$ is not congruent to $-C^2 - \frac{3}{2}(1 + b_+(X))$ modulo 4, and set $\Phi_C^X := \sum_{d \geq 0} \Phi_{C,d}^X$, thought of as a linear map from $A_*(X) = \sum_{d \geq 0} A_d(X)$ to $\mathbb{Q}$. The polynomials $\Phi_{C,d}^X$ depend only on $C$ modulo $4H^2(X, \mathbb{Z})$ and up to sign only on $C$ modulo $2H^2(X, \mathbb{Z})$.

Kronheimer and Mrowka introduce the notion of *simple type*: The manifold $X$ is of simple type if $\Phi_C^X(\alpha(p^2 - 4)) = 0$ for all $\alpha \in A_*(X)$ and all $C \in H^2(X, \mathbb{Z})$. In [Kr-M1], [Kr-M2], [F-S2] the structure of the Donaldson invariants of 4-manifolds of simple type is then described in a very compact way: There exist finitely many so-called basic classes $K_i \in H^2(X, \mathbb{Z})$ (independent of $C$) and rational numbers $a_i$ (depending on $C$) such that for all $x \in H_2(X, \mathbb{Z})$ there is an identity

$$\Phi_C^X\big((1 + p/2)e^{xz}\big) = e^{\mathbf{Q}(x)z^2/2} \sum_i a_i e^{K_i x z},$$

of formal power series in $z$, where $\mathbf{Q}$ is the quadratic form on $H_2(X, \mathbb{Z})$. Note that this formula determines $\Phi_C^X$ completely since $\Phi_{C,d}^X(p^i x^n)$ is 0 unless $2i + n \equiv -C^2 - \frac{3}{2}(1 + b_+(X))$ modulo 4. Very many 4-manifolds with $b_+ > 1$ have been shown to be of simple type, and the class of 4-manifolds of simple type is closed under several natural operations on 4-manifolds, like connected sum with $\bar{\mathbb{P}}_2$ and rational blowdown [F-S3].

Fintushel and Stern and Kronheimer and Mrowka have also introduced a generalization of the concept of simple type, the *k-th order simple type*, where one instead requires that $\Phi_C^X((p^2 - 4)^k \alpha) = 0$. They showed that all manifolds with $b_+ > 1$ are of $k$-th order simple type for some $k$, but it is not known whether $k > 1$ is ever needed.

A new light has been shed on these structure theorems by the Seiberg-Witten invariants [S-W],[W1]. A class $K$ in $H^2(X, \mathbb{Z})$ is called *SW-basic* if the corresponding Seiberg-Witten invariant does not vanish, and $X$ is *of SW-simple type* if for all $SW$-basic classes the corresponding moduli space is 0-dimensional. It is conjectured that the $SW$-basic classes are the same as the basic classes in Donaldson theory, and that the condition of simple type in Donaldson theory and in Seiberg-Witten theory are equivalent. From the viewpoint of theoretical physics the precise relation between the Donaldson and the Seiberg-Witten invariants should be given in terms of the modular curve $\mathfrak{H}/\Gamma(2)$.

There is a project towards giving a mathematical proof of this relationship [P-T1], [P-T2]. All symplectic 4-manifolds with $b_+ > 1$ are known to be of $SW$-simple type [T2], and no 4-manifold with $b_+ > 1$ is known not to be of $SW$-simple type.

In this paper we study the case $b_+ = 1$, where the invariants have been introduced in [K1]. In this case the invariants are no longer independent of the metric, but in [K-M] they were shown to depend only on the so-called *period point* of the metric in the positive cone $H^2(X,\mathbb{R})^+$. We will therefore denote them by $\Phi_{C,d}^{X,L}$, where $L$ is the period point. In fact there is a collection of cohomology classes $\xi \in H^2(X,\mathbb{Z})$ defining "walls" (hyperplanes), and $\Phi_{C,d}^{X,L}$ depends only on the chamber (= connected component of the complement of the walls) to which $L$ belongs. The change of the Donaldson invariants when passing through a wall is given by *wall-crossing terms* $\delta_{\xi,d}^X$. Kotschick and Morgan made a conjecture about the structure of the $\delta_{\xi,d}^X$. In particular they should only depend on $\xi^2$, $d$ and the homotopy type of $X$. Morgan and Ozsváth have announced a proof of this conjecture [M-O]. In this whole paper we will assume the conjecture.

In [K-L] the wall-crossing formulas were used in combination with the blowup formulas to compute Donaldson invariants of $\mathbb{P}_2$ and $\mathbb{P}_1 \times \mathbb{P}_1$ and to show in particular that $\mathbb{P}_2$ is not of simple type and $\mathbb{P}_1 \times \mathbb{P}_1$ is not of simple type for any chamber. In [E-G1],[F-Q],[E-G2] the $\delta_{\xi,d}^S$ were studied in the case of algebraic surfaces $S$ with $p_g = 0$. In [E-G1],[F-Q] they were expressed in terms of Hilbert schemes of points on $S$, and the leading terms were computed. In [E-G2] the Bott residue formula is applied to computing the $\delta_{\xi,d}^S$, and the Donaldson invariants, for rational surfaces with help of the computer. In [G], assuming the conjecture of Kotschick and Morgan and using the blowup formulas, the $\delta_{\xi,d}^X$ were determined completely (for arbitrary $X$) in terms of modular forms. On a rational algebraic surface one can, after possibly blowing up, always find a chamber where the Donaldson invariants vanish. Therefore the Donaldson invariants of rational surfaces can always be expressed in terms of modular forms. In particular this gives the Donaldson invariants of $\mathbb{P}_2$. The results made it look very unlikely that $\mathbb{P}_2$ or $\mathbb{P}_1 \times \mathbb{P}_1$ can be of $k$-th order simple type for any $k$. By just straightforwardly applying the results from [G] one will in general get a very complicated expression for the Donaldson invariants, from which it is very hard to read off structural results.

On the other hand Morgan and Szabó showed [M-Sz] that, for some rational surfaces admitting an elliptic fibration with $CF$ odd ($F \in H^2(X,\mathbb{Z})$ dual to the class of a fibre), the *limit* of the Donaldson invariants $\Phi_C^{X,L}$, for the period point $L$ tending to $F$, fulfills the simple type condition $\Phi_C^{X,L}(\alpha(p^2 - 4)) = 0$, and the corresponding Donaldson power series is given by formulas analoguous to those of [Kr-M2] and [F-S2]. In fact this is a special case of their more general results.

In the current paper we want to apply the results of [G] to understand the structure of the Donaldson invariants of 4-manifolds with $b_+ = 1$. Looking at the previously known results (e.g. [K-L],[E-G2]), it seems very unlikely that results analogous to [Kr-M2] and [F-S2] can hold for period points in the *inside* of the positive cone, while the results of [M-Sz] and the structure of the formulas in [G] suggest that one should restrict the attention to the *boundary*. One first has to find the correct definition of the Donaldson invariants for a period point $F$ there. If, for a primitive representative $F \in H^2(X,\mathbb{Z})$, the number $CF$ is odd, then $F$ lies in the closure of a unique chamber and on no wall, and one just takes the values of the invariants in this chamber. In case $CF$ is even, $F$ will in general lie in the closure of infinitely many chambers, and one has to take a certain kind of "renormalized average" over all of them.

With this correction the Donaldson invariants are well-defined, and we show that for $F, G \in H^2(X, \mathbb{Z})$ on the boundary of the positive cone, the *difference* $\Phi_C^{X,F} - \Phi_C^{X,G}$ always fullfils the $k$-th order simple type condition for $k = (W^2 - \sigma(X))/8$, where $\sigma(X)$ is the signature of $X$ and $W$ is a characteristic element with $W^2$ maximal in a certain sector of $H^2(X, \mathbb{Z})$ determined by $F$ and $G$. We also give a universal formula for the precise structure of this difference in terms of modular forms and an explicit set of *basic* cohomology classes, which are again characteristic elements in the above sector of $H^2(X, \mathbb{Z})$. With $k$ as above the leading terms of this formula give an expression for $(\Phi_C^{X,F} - \Phi_C^{X,G})(e^{xz}(1+p/2)(1-p^2/4)^{k-1})$ analogous to that of [Kr-M2] and [F-S2]. There is however one difference: While in the case $b_+ > 1$ there is only a finite number of basic classes, one can interpret the formulas in our case as saying that there are infinitely many (all but finitely many orthogonal to $F$ or $G$), and that the Donaldson invariants are obtained as a "renormalized" (by analytic continuation) sum of their contributions.

We note that the basic classes not orthogonal to $F$ or $G$ are precisely the characteristic cohomology classes whose corresponding Seiberg-Witten invariants (with respect to the unperturbed equations) for metrics with period points near $F$ and near $G$ differ. This leads to a conjectural formula for the relationship between Seiberg-Witten and Donaldson invariants in the case $b_+ = 1$ on the boundary of the positive cone.

If $X$ is a rational algebraic surface, it is easy to show that there always exists a point on the boundary of the postive cone such that all Donaldson invariants vanish. Therefore for these surfaces the formulas above give the Donaldson invariants and not just their differences. In particular the above-mentioned conjecture about the relation of Donaldson and Seiberg-Witten invariants holds for rational surfaces. In particular this shows that for any $k$ very many $(X, F)$ are of strictly $k$-th order simple type, related to the fact that in case $b_+ = 1$, the chamber structure in Seiberg-Witten theory will very often make $SW$-basic classes of arbitrary expected dimension appear.

In the case of $\mathbb{P}_2 \# \overline{\mathbb{P}}_2$ or $\mathbb{P}_1 \times \mathbb{P}_1$, the limit of the Donaldson invariants at the boundary of the positive cone is also defined without the correction, but the results of [K-L] show that they are not of simple type, whereas with the correction they vanish. The limit can therefore be computed as the negative of the correction we introduce, proving a conjecture from [E-G2]. We conclude by giving a number of examples for our formulas.

To prove our results we observe that by the results of [G] the formula for the difference of the Donaldson invariants at two boundary points is very closely related to a new kind of theta functions, associated to a lattice $L$ of type $(r-1, 1)$ and a pair $(f, g)$ of elements in $L$ with self-intersection 0 (in fact the lattice is $H_2(X, \mathbb{Z})$ with the negative of the intersection form). We show that these theta functions are Jacobi forms for the theta group $\Gamma_\theta$. Using this fact one can show that the difference of the Donaldson invariants has a development in terms of modular functions for a subgroup $\Gamma_u \subset \Gamma_\theta$ which are holomorphic except at the cusps $-1$, $1$, $\infty$ of $\Gamma_u$. It follows that they can be expressed as a rational function in a certain modular function $u$, and that the relevant information about the Donaldson invariants can be read off from the poles at the cusp $-1$. We note that $\Gamma_u$ is conjugate to $\Gamma(2)$ by an element of $GL(2,\mathbb{Z})$, which maps the cusps $-1, 1, \infty$ to $\infty, 0, 1$, thus giving a connection to the description from theoretical physics. In particular all the computations could be rephrased in terms of $\Gamma(2)$.

We also mention the connection with the recent work of Borcherds on automorphic forms on orthogonal groups ([Bo], in particular the results of §10 and the examples and problems concerning Donaldson polynomials in §15 and §16).

The first named author would like to thank Barbara Fantechi, with whom he had useful discussions on several aspects of this paper, Ronald Stern, who told him about the notion of higher order simple type, Victor Pidstrigach and Andrei Tyurin who explained their program for proving the relationship of Donaldson and Seiberg-Witten invariants to him, and John Morgan and Zoltan Szabó, who informed him about their results about the structure of some Donaldson invariants e.g. on $\mathbb{P}_2 \# 9\overline{\mathbb{P}}_2$; in fact this paper was motivated in part by trying to understand that result. With Zoltan Szabó he also had some further discussions, which helped to clarify our ideas. This work was started while the first named author was at the Max-Planck-Institut für Mathematik Bonn, and carried out during his stay at Pisa, with a grant of MAP.

## 2. Preliminaries

In this paper let $X$ be a simply connected smooth 4-manifold with $b_+(X) = 1$.

**Notation 2.1.** We will usually denote by upper case letters the classes in $H^2(X,\mathbb{C})$, unless these appear as walls (see below), when we denote them by Greek letters. Elements of $H_2(X,\mathbb{C})$ are denoted by lower case letters. We usually denote the Poincaré dual class of $A \in H^2(X,\mathbb{C})$ by the corresponding lower case letter $a$. For $A, B \in H^2(X,\mathbb{C})$ and the corresponding dual classes $a, b \in H_2(X,\mathbb{C})$, the canonical pairing of $H^2(X,\mathbb{C})$ and $H_2(X,\mathbb{C})$ and the intersection product on $H^2(X,\mathbb{C})$ are just denoted by $Ab$ and $AB$ respectively. We denote $\mathbf{Q}(a)$ the quadratic form (given by $A^2$) on $H_2(X,\mathbb{C})$ and by $\sigma(X)$ the signature of $X$. We denote $\widehat{X} := X \# \overline{\mathbb{P}}_2$, and $E$ the canonical generator of $H^2(\overline{\mathbb{P}}_2, \mathbb{Z})$. We will always identify $H^2(X,\mathbb{Q})$ with the orthogonal subspace $E^\perp \subset H^2(\widehat{X},\mathbb{Q})$. We trust that there will be no confusion between the exponential function and the Poincaré dual of $E$.

Let $\mathcal{E}$ a complex rank 2 bundle with first Chern class $C$ and second Chern class $c_2$. We put $d := 4c_2 - C^2 - 3$, and denote by $\Phi_{C,d}^{X,g}$ the Donaldson invariant corresponding to $\mathcal{E}$ and the (generic) metric $g$ (cf. [Do]). They depend up to sign only on $C$ modulo 2, and the sign is determined by $C$ modulo 4 and by a so-called *homology orientation*. We use the conventions of e.g. [F-S1]. Let $p \in H_0(X,\mathbb{Z})$ be the class of a point. Let $A_d(X)$ be the set of polynomials of weight $d$ in $H_2(X,\mathbb{C}) \oplus H_0(X,\mathbb{C})$, where $a \in H_2(X,\mathbb{C})$ has weight 1 and $p$ has weight 2, and put $A_*(X) := \bigoplus_{d \geq 0} A_d(X)$. Then $\Phi_{C,d}^{X,g}$ is a linear map $A_d(X) \longrightarrow \mathbb{Q}$. The letters $z$ and $t$ will denote indeterminates; $z$ will also be sometimes a complex variable, and $x$ usually denotes an element of $H_2(X,\mathbb{C})$. We will very often write $1^a$ instead of $e^{2a\pi i}$ for a rational number $a$.

2.1. **Walls and chambers.** In the case $b_+ = 1$ the Donaldson invariants are no longer independent of the metric.

**Definition 2.2.** Firstly we recall that a Riemannian metric $g$ determines a ray in
$$H^2(X,\mathbb{R})^+ = \{H \in H^2(X,\mathbb{R}) \mid H^2 > 0\};$$
namely the set of self-dual harmonic forms. This ray, or any representative in $H^2(X,\mathbb{R})^+$, is called the *period point* of $g$ and denoted $\omega(g)$. The quotient $H^2(X,\mathbb{R})^+/\mathbb{R}^+$ has two connected components. The choice of a homology orientation amounts to the choice of one of them, which we call $\mathbb{H}_X$. We will in future always assume that we have chosen a homology orientation, and view the period points as lying in the corresponding component $\mathbb{H}_X$.

The space $\mathbb{H}_X$ is a model of hyperbolic $r$-space, where $r = b_2(X) - 1$. In particular we can complete it to $\overline{\mathbb{H}}_X := \mathbb{H}_X \cup \mathbb{S}_X$, where
$$\mathbb{S}_X := (\{H \in H^2(X,\mathbb{Q}) \mid H^2 = 0\} \setminus \{0\})/\mathbb{Q}^*$$
is the set of *cusps*. We will usually not distinguish between an element $H \in \mathbb{H}_X$ and a representative in $H^2(X,\mathbb{R})^+$, and similar for $\mathbb{S}_X$ (with the representative on the boundary of the chosen component).

**Definition 2.3.** (see e.g. [K1], [K-M]) By a *wall* in $\mathbb{H}_X$ we mean the intersection of $\mathbb{H}_X$ with a set
$$W^\xi := \{L \in H^2(X,\mathbb{R})^+ \mid \xi L = 0\}/\mathbb{R}^+.$$
where $\xi \in H^2(X,\mathbb{Q})$ with $\xi^2 < 0$. If $C \in H^2(X,\mathbb{Z})$, and $d \in \mathbb{Z}_{\geq 0}$, we will say that an element $\xi \in H^2(X,\mathbb{Q})$ is *of type* $(C,d)$ if $C/2 - \xi \in H^2(X,\mathbb{Z})$ and $(d+3)/4 + \xi^2 \in \mathbb{Z}_{\geq 0}$. A *chamber* of type $(C,d)$ is a connected component of the complement of the walls in $\mathbb{H}_X$ defined by elements of type $(C,d)$.

**Theorem 2.4.** [K-M]
1. $\Phi_{C,d}^{X,g}$ depends only on the chamber of the period point $\omega(g)$.
2. For all $\xi \in \frac{1}{2}H^2(X,\mathbb{Z})$ of type $(C,d)$ there exists a linear map $\delta_{\xi,d}^X : A_d(X) \longrightarrow \mathbb{C}$ such that for all generic metrics $g_1$ and $g_2$
$$\Phi_{C,d}^{X,g_1} - \Phi_{C,d}^{X,g_2} = 1^{C^2/8} \sum_\xi (-1)^{(\xi - C/2)C} \delta_{\xi,d}^X,$$
where the sum runs through all $\xi$ of type $(C,d)$ with $\xi\omega(g_2) < 0 < \xi\omega(g_1)$.

Note that the conventions are different from [G],[K-M],[K-L]. We have changed the sign conventions and replaced $\xi$ by $\xi/2$ to get below a more direct relation to theta functions.

Actually in [K-M] the result is only proven for the restriction of $\Phi_{C,d}^{X,g}$ to $\operatorname{Sym}^d(H_2(X,\mathbb{Q}))$; see [G] for the extension to $A_d(X)$. We mention that this extension also works outside of the so-called *stable range*, i.e. when $C \equiv 0$ modulo $2$ and for classes $\alpha \in A_d(X)$ containing monomials $x^k p^r$ with $x \in H_2(X,\mathbb{Q})$ and $k < 2(r+1)$.

Kotschick and Morgan make a conjecture about the structure of the wallcrossing terms.

**Conjecture 2.5.** [K-M] $\delta^X_{\xi,d}|_{Sym^d(H_2(X,\mathbb{Q}))}$ *is a polynomial in $\xi$ and the quadratic form $\mathbf{Q}$ whose coefficients depend only on $\xi^2$, $d$, and the homotopy type of $X$.*

In this statement, polynomials in $\xi$ and $\mathbf{Q}$ are considered as maps on $\operatorname{Sym}^d(H_2(X),\mathbb{Q})$ by $\mathbf{Q}^k\xi^\ell(x^d) = \mathbf{Q}(x)^k(\xi x)^\ell$ if $2k+\ell = d$ (and 0 otherwise) for $x \in H_2(X,\mathbb{Q})$, and then by multilinearity for arbitrary elements of $\operatorname{Sym}^d(H_2(X),\mathbb{Q})$.

A proof of Conjecture 2.5 has been announced [M-O]. We will assume 2.5 for the whole of the paper.

By Theorem 2.4 the Donaldson invariant $\Phi^{X,g}_{C,d}$ depends only on the chamber of the period point $\omega(g) \in \mathbb{H}_X$. For $H \in \mathbb{H}_X$ not on a wall defined by a class of type $(C,d)$ we will therefore put $\Phi^{X,H}_{C,d} = \Phi^{X,g}_{C,d}$, where $g$ is a generic metric whose period point lies in the same chamber of type $(C,d)$ as $H$. Apparently it is not known whether every element of $\mathbb{H}_X$ appears as period point of a metric. If there is no period point in the chamber of $H$, then we define $\Phi^{X,H}_{C,d}$ by requiring that

$$\Phi^{X,H}_{C,d} - \Phi^{X,g}_{C,d} = 1^{C^2/8} \sum_\xi (-1)^{(\xi - C/2)C} \delta^X_{\xi,d}, \tag{2.5.1}$$

the sum running through all $\xi$ of type $(C,d)$ with $\xi\omega(g) < 0 < \xi H$, where $g$ is any generic Riemannian metric on $X$. Theorem 2.4 implies that in this way $\Phi^{X,H}_{C,d}$ is well-defined.

Finally, we mention that if $C \equiv 0$ modulo 2 then $\Phi^{X,H}_{0,d}(\alpha)$ was originally defined only for $\alpha$ in the stable range but can be defined for arbitrary $\alpha$ by

$$\Phi^{X,H}_{C,d}(\alpha) := \Phi^{\widehat{X},H+\epsilon E}_{C+E,d+1}(e\alpha) \qquad (\epsilon > 0 \text{ sufficently small}) \tag{2.5.2}$$

(see [G] definition 2.7). Note that now by the blowup formulas (e.g. [F-S1],[T1]) the formula (2.5.2) holds for all $d$ and all $\alpha \in A_d(X)$.

We put $\Phi^{X,H}_{C,d} = 0$ if $d$ is not congruent to $-C^2 - 3$ modulo 4 and, if $H\xi \neq 0$ for all $\xi \in H^2(X,\mathbb{Z}) + C/2$ (i.e. $H$ does not lie on a wall defined by a class of type $(C,d)$ for any $d$), we put $\Phi^{X,H}_C := \sum_{d\geq 0} \Phi^{X,H}_{C,d}$. We also put $\delta^X_\xi := \sum_{d\geq 0} \delta^X_{\xi,d}$. For an indeterminate $z$ and $a_n \in A_*(X)$ for all $n \geq 0$ we put

$$\Phi^{X,H}_C\left(\sum_{n\geq 0} a_n z^n\right) := \sum_{n\geq 0} \Phi^{X,H}_C(a_n) z^n, \quad \delta^X_\xi\left(\sum_{n\geq 0} a_n z^n\right) := \sum_{n\geq 0} \delta^X_\xi(a_n) z^n$$

With this convention we define, for $x \in H_2(X,\mathbb{Z})$ and $P(p)$ a polynomial in $p$ (the class of a point) the following two formal power series in variables $z$, $t$:

$$\Psi^{X,H}_C(x \cdot z, P(p)) := \Phi^{X,H}_C(e^{xz} P(p)), \quad \overline{\Psi}^{X,H}_C(x \cdot z, t) := \sum_{r=0}^\infty \Psi^{X,H}_C(x \cdot z, p^r) t^{r+1}. \tag{2.5.3}$$

(The change of argument from $e^{xz}$ in $\Phi$ to $x \cdot z$ in $\Psi$ will become important in Section 4, where the invariant $\Psi$ —but not $\Phi$—is defined for $H$ belonging to the boundary of the positive cone.)

Following [Kr-M1],[Kr-M2] we say that the pair $(X, H)$ is *of $k$-th order simple type* if $\Phi^{X,H}_C(\alpha(p^2 - 4)^k) = 0$ for all $\alpha \in A_*(X)$ and all $C \in H^2(X,\mathbb{Z})$. It is of strictly $k$-th order simple type if it is of $k$-th order simple type but not of $(k-1)$-th order simple type.

**2.2. Some elementary computations with modular forms.** Let $\mathfrak{H} := \{\tau \in \mathbb{C} \mid \mathrm{Im}(\tau) > 0\}$ be the complex upper half-plane. For $\tau \in \mathfrak{H}$ let $q := e^{2\pi i \tau}$ and $q^{1/n} := e^{2\pi i \tau/n}$. Elements $A = \begin{pmatrix} a & b \\ c & d \end{pmatrix} \in SL(2, \mathbb{Z})$ act on $\mathfrak{H}$ by $A\tau := \dfrac{a\tau + b}{c\tau + d}$ and, given $k \in \mathbb{Z}$, on functions $g : \mathfrak{H} \longrightarrow \mathbb{C}$ by

$$(g|_k A)(\tau) := (c\tau + d)^{-k} g(A\tau).$$

For $k \in \tfrac{1}{2}\mathbb{Z}$ and $(c,d) \neq (0,-1)$ we define $(g|_k A)$ by the same formula, where $(c\tau+d)^{-k}$ stands for $\sqrt{c\tau+d}^{-2k}$, where we always use the principal branch of the square root (whose real part is positive on complex numbers with argument strickly between $-\pi$ and $\pi$). Note that this in general does not define an action anymore. If $g$ is a modular form of weight $k$ or $k$ is otherwise understood we just write $g|A$. Let

$$T := \begin{pmatrix} 1 & 1 \\ 0 & 1 \end{pmatrix}, \quad V := T^2 = \begin{pmatrix} 1 & 2 \\ 0 & 1 \end{pmatrix}, \quad S := \begin{pmatrix} 0 & -1 \\ 1 & 0 \end{pmatrix}, \quad W := T^{-1}S = \begin{pmatrix} -1 & -1 \\ 1 & 0 \end{pmatrix}. \tag{2.5.4}$$

Let $\Gamma_\theta = \langle V, S \rangle \subset SL(2, \mathbb{Z})$ be the theta group. The quotient $\mathfrak{H}/\Gamma_\theta$ has two cusps: $-1$ and $\infty$. Let $\eta(\tau) := q^{1/24} \prod_{n>0}(1-q^n)$ be the Dirichlet eta function and $\Delta := \eta^{24}$ the discriminant. We denote $\sigma_k(n) := \sum_{d|n} d^k$ and by $\sigma_1^{\mathrm{odd}}(n)$ the sum of the odd divisors of $n$. For even $k \geq 2$ let

$$G_k(\tau) := -\frac{B_k}{2k} + \sum_{n>0} \sigma_{k-1}(n) q^n$$

be the classical Eisenstein series, where $B_k$ is the $k$-th Bernoulli number. For odd $k$ we put $G_k(\tau) := 0$. Note that $G_k$ is a modular form of weight $k$ for $SL(2, \mathbb{Z})$ for $k \geq 4$, but is only "quasi-modular" for $k = 2$, i.e., it transforms by equation 3.11.1.

Recall the classical theta functions

$$\theta_{\mu\nu}(\tau, z) := \sum_{n \in \mathbb{Z}} (-1)^{n\nu} q^{(n+\mu/2)^2/2} e^{2\pi i (n+\mu/2) z} \quad (\mu, \nu \in \{0,1\}) \tag{2.5.5}$$

and their "Nullwerte"

$$\theta(\tau) := \theta_{00}(\tau, 0) = \frac{\eta(\tau)^5}{\eta(\tau/2)^2 \eta(2\tau)^2}, \quad \theta_{01}^0(\tau) = \theta_{01}(\tau, 0) = \frac{\eta(\tau/2)^2}{\eta(\tau)}$$

$$\theta_{10}^0(\tau) = \theta_{10}(\tau, 0) = 2\frac{\eta(2\tau)^2}{\eta(\tau)}, \quad \theta_{11}^0(\tau) = 0.$$

We set

$$f(\tau) := 1^{-1/8} \frac{\eta(\tau)^3}{\theta(\tau)} = 1^{-1/8} \frac{\eta(\tau/2)^2 \eta(2\tau)^2}{\eta(\tau)^2} = \frac{1^{-1/8}}{2} \theta_{01}^0(\tau) \theta_{10}^0(\tau). \tag{2.5.6}$$

The transformation laws

$$\eta(\tau + 1) = 1^{1/24} \eta(\tau), \qquad \eta(-1/\tau) = \sqrt{\frac{\tau}{i}} \eta(\tau)$$

(where we again use the principal branch of the square root) imply the following.

1. $\theta|V = \theta$, $\theta|S = 1^{-1/8}\theta$, $\theta|T = \theta_{01}^0$ and $\theta|W(\tau) = 1^{-1/8}\theta_{10}^0(\tau) = 1^{-1/8} \frac{2\eta(2\tau)^2}{\eta(\tau)}$.
2. $f|V = if$, $f|S = -if$, $f|T^{-1} = 1^{-1/4} \frac{\eta(\tau)^3}{\theta_{01}^0}$ and $f|W(\tau) = -\frac{\eta(\tau)^4}{2\eta(2\tau)^2}$. In particular $f^4$ is a modular form of weight 4 for $\Gamma_\theta$.

3. The function
$$R(\tau) := \frac{\Delta(\tau)^2}{\Delta(\tau/2)\Delta(2\tau)}$$
is a modular function for $\Gamma_\theta$ and $R|W(\tau) = -2^{12}\Delta(2\tau)/\Delta(\tau)$.

Finally let $e_1, e_2, e_3$ be the 2-division values of the Weierstraß $\wp$-function at $1/2$, $\tau/2$ and $(1+\tau)/2$ respectively, i.e.

$$e_1(\tau) = -\frac{1}{6} - 4\sum_{n>0}\sigma_1^{\mathrm{odd}}(n)q^n,$$

$$e_2(\tau) = \frac{1}{12} + 2\sum_{n>0}\sigma_1^{\mathrm{odd}}(n)q^{n/2},$$

$$e_3(\tau) = \frac{1}{12} + 2\sum_{n>0}(-1)^n\sigma_1^{\mathrm{odd}}(n)q^{n/2}.$$

We write
$$U(\tau) := -\frac{3e_3(\tau)}{f(\tau)^2}, \quad u(\tau) := \frac{1}{U(\tau)}, \quad G(\tau) := G_2(\tau) + e_3(\tau)/2. \tag{2.5.7}$$

The first two are modular functions, the third a quasi-modular form.

**Lemma 2.6.** 1. *$e_3$ is a modular form of weight 2 for $\Gamma_\theta$ and $e_3|T = e_2$, $e_3|W = e_1$.*

2. *We have the identities*

$$U(\tau)^2 - 4 = -\frac{1}{16}R(\tau), \tag{2.6.1}$$

$$q\frac{d}{dq}U(\tau) = -\frac{1}{16}R(\tau)f(\tau)^2. \tag{2.6.2}$$

*Proof.* 1. is standard and can be shown e.g. by observing that $e_3(\tau) = -2G_2(\tau/2) + 8G_2(\tau) - 8G_2(2\tau)$ and using the transformation behaviour of the quasi-modular form $G_2$. For 2., we note that the four functions $f^4$, $e_3^2 = U^2 f^4/9$, $Rf^4$ and $f^2 q\frac{d}{dq}U$ are modular forms of weight 4 on $\Gamma_\theta$. Therefore they are linear combinations of the two Eisenstein series $G_4(\tau)$ and $G_4((\tau+1)/2)$, which generate the space of modular forms of weight 4 for $\Gamma_\theta$. The identities (2.6.1) and (2.6.2) are then obtained by comparing the first few Fourier coefficients. □

**Remark 2.7.** Let $\Gamma_u = \pm\langle V^2, VS, SV\rangle$; this is a subgroup of index 2 of $\Gamma_\theta$. The quotient $\mathfrak{H}/\Gamma_u$ has 3 cusps $1$, $-1$ and $\infty$. From the above it is clear that $U(\tau)$ is a modular function for $\Gamma_u$ and that $f(\tau)^2$ is a modular form of weight 2 for $\Gamma_u$. We note that $U(\tau)$ defines an isomorphism of from $\mathfrak{H}/\Gamma_u \cup \{1\} \cup \{-1\} \cup \{\infty\}$ to $\mathbb{P}_1$, as the smallest power of $T$ contained in $\Gamma_u$ is $T^4$ and the lowest power of $q$ in the Fourier development of $U(\tau)$ is $-1/4$. We see that $U(\infty) = \infty$, $U(-1) = 2$, $U(1) = -2$. In fact $W$ maps $\infty$ to $-1$ and $VW$ maps $\infty$ to $1$, and

$$\widetilde{U}(\tau) := U|W(\tau) = -\frac{12e_1(\tau)\eta(2\tau)^4}{\eta(\tau)^8} = 2 + 64q + 512q^2 + 2816q^3 + \ldots, \quad U|VW(\tau) = -\widetilde{U}(\tau).$$

We denote $\widetilde{u}(\tau) = u|W(\tau) = 1/\widetilde{U}(\tau)$. In particular every modular function for $\Gamma_u$ will be a rational function in $U(\tau)$, and a polynomial in $U(\tau)$, $(U(\tau) - 2)^{-1}$ and $(U(\tau) + 2)^{-1}$ if it is holomorphic on $\mathfrak{H}$. We note that $\Gamma_u$ is conjugate to $\Gamma(2)$ via the matrix $\begin{pmatrix} 1 & -1 \\ 1 & 1 \end{pmatrix} \in GL(2, \mathbb{Q})$, which sends the cusps $\infty, 1, -1$ to $1, 0, \infty$ respectively.

We finish this section by giving a list of the leading terms of the Fourier developments of some of the modular forms and functions we introduced, which will play a role in our later computations. This should help the reader to check and apply our computations and results.

$$\theta(\tau) = 1 + 2q^{\frac{1}{2}} + 2q^2 + 2q^{\frac{9}{2}} + 2q^8 + \ldots$$

$$f(\tau) = 1^{-\frac{1}{8}}q^{\frac{1}{8}}(1 - 2q^{\frac{1}{2}} + q - 2q^{\frac{3}{2}} + 2q^2 + 3q^3 - 2q^{\frac{7}{2}} - 2q^{\frac{9}{2}} + 2q^5 - 2q^{\frac{11}{2}} + q^6 - 2q^{\frac{13}{2}} + \ldots)$$

$$R(\tau) = q^{-\frac{1}{2}} + 24 + 276q^{\frac{1}{2}} + 2048q + 11202q^{\frac{3}{2}} + 49152q^2 + 184024q^{\frac{5}{2}} + 614400q^3$$
$$+ 1881471q^{\frac{7}{2}} + 5373952q^4 + 14478180q^{\frac{9}{2}} + \ldots$$

$$e_1(\tau) = -\frac{1}{6} - 4q - 4q^2 - 16q^3 - 4q^4 - 24q^5 - 16q^6 - 32q^7 - 4q^8 - \ldots$$

$$e_3(\tau) = \frac{1}{12} - 2q^{\frac{1}{2}} + 2q - 8q^{\frac{3}{2}} + 2q^2 - 12q^{\frac{5}{2}} + 8q^3 - 16q^{\frac{7}{2}} + 2q^4 - 26q^{\frac{9}{2}} + 12q^5 - 24q^{\frac{11}{2}} + \ldots$$

$$U(\tau) = 1^{\frac{1}{4}}q^{-\frac{1}{4}}(-\frac{1}{4} + 5q^{\frac{1}{2}} + \frac{31}{2}q + 54q^{\frac{3}{2}} + \frac{641}{4}q^2 + 409q^{\frac{5}{2}} + \frac{1889}{2}q^3 + 2062q^{\frac{7}{2}}$$
$$+ \frac{17277}{4}q^4 + 8666q^{\frac{9}{2}} + \frac{33439}{2}q^5 + 31328q^{\frac{11}{2}} + 57313q^6 + \ldots)$$

$$u(\tau) = 1^{\frac{1}{4}}q^{\frac{1}{4}}(4 + 80q^{\frac{1}{2}} + 1848q + 42784q^{\frac{3}{2}} + 990100q^2 + 22911600q^{\frac{5}{2}} + 530190104q^3$$
$$+ 12268965984q^{\frac{7}{2}} + 283912371144q^4 + \ldots)$$

$$G(\tau) = -q^{\frac{1}{2}} + 2q - 4q^{\frac{3}{2}} + 4q^2 - 6q^{\frac{5}{2}} + 8q^3 - 8q^{\frac{7}{2}} + 8q^4 - 13q^{\frac{9}{2}} + 12q^5 - 12q^{\frac{11}{2}} + 16q^6 + \ldots$$

$$\widetilde{U}(\tau) = 2 + 64q + 512q^2 + 2816q^3 + 12288q^4 + 45952q^5 + 153600q^6 + 470528q^7$$
$$+ 1343488q^8 + 3619136q^9 + 9280512q^{10} + \ldots$$

$$\theta_{10}^0(\tau) = q^{\frac{1}{8}}(2 + 2q + 2q^3 + 2q^6 + 2q^{10} + \ldots)$$

$$\frac{\eta(\tau)^4}{\eta(2\tau)^2} = 1 - 4q + 4q^2 + 4q^4 - 8q^5 + 4q^8 - 4q^9 + 8q^{10} + \ldots$$

2.3. **Wall-crossing formula.** In [G] definition 2.7 we extended the definition of the wall-crossing terms $\delta_{\xi,d}^X$ to all classes $\xi \in H^2(X, \mathbb{Z})$ by use of the blowup formulas [F-S1],[T1]. If $\xi^2 < 0$ we have in particular $\delta_{\xi,d}^X = 0$, if $\xi$ is not of of type $(2\xi, d)$. The main theorem of [G] is:

**Theorem 2.8.** [G] *Let $X$ be a simply connected 4-manifold with $b_+ = 1$ and signature $\sigma(X)$. For $x \in H_2(X, \mathbb{Z})$ put*

$$\Delta_\xi^X(\tau, x \cdot z) := -\frac{4\theta(\tau)^{\sigma(X)}}{f(\tau)} \cdot q^{-\xi^2/2} e^{-\xi x z/f(\tau)} \cdot e^{-\mathbf{Q}(x)G(\tau)z^2/f(\tau)^2}.$$

*Then*

$$\delta_\xi^X(e^{xz}p^r) = \text{Coeff}_{u(\tau)^{r+1}}\left[\Delta_\xi^X(\tau, x \cdot z)\right]. \tag{2.8.1}$$

Here the symbol "Coeff" means that we expand the following expression (or more precisely, the coefficient of $z^n$ in it for each fixed power $n$) as a fractional Laurent series in $q$, then rewrite it as as a fractional Laurent series in $u(\tau) = 4iq^{\frac{1}{4}} + \cdots$, and then take the coefficient of the indicated power of $u(\tau)$.

Theorem 2.8 was stated in [G] (up to some differences in conventions) in the form

$$\delta_\xi^X \left( e^{xz} p^r \right) := \operatorname{Res}_{q=0} \left[ \frac{1}{4} f(\tau)^2 U(\tau)^r R(\tau) \Delta_\xi^X(\tau, x \cdot z) \frac{dq}{q} \right]. \tag{2.8.2}$$

Equation (2.8.1) follows from this by

**Proposition 2.9.** *Let $F$ be a meromorphic function on $\mathfrak{H}$ having a Laurent development in powers of $q^{1/n}$ for some $n \in \mathbb{Z}_{>0}$. Then*

$$\operatorname{Res}_{q=0} \left[ \frac{1}{4} f(\tau)^2 U(\tau)^r R(\tau) F(\tau) \frac{dq}{q} \right] = \operatorname{Coeff}_{u(\tau)^{r+1}} \left[ F(\tau) \right].$$

*Proof.* We use the formula (2.6.2) to obtain

$$\frac{1}{4} \operatorname{Res}_{q=0} \left[ f(\tau)^2 \frac{R(\tau)}{u(\tau)^r} F(\tau) \frac{dq}{q} \right] = 4 \operatorname{Res}_{q=0} \left[ F(\tau) \left( \frac{d}{dq} u(\tau) \right) / u(\tau)^{r+2} \, dq \right]$$

$$= \operatorname{Res}_{u(\tau)=0} \left[ F(\tau) \frac{d\,u(\tau)}{u(\tau)^{r+2}} \right] = \operatorname{Coeff}_{u(\tau)^{r+1}} \left[ F(\tau) \right].$$

$\square$

**Remark 2.10.** 1. The term $q^{-\xi^2/2} e^{-\xi xz}$ occuring in $\Delta_\xi^X(\tau, x \cdot z)$ looks like the summand corresponding to a vector $\xi \in L$ in the theta series for the lattice $L = H_2(X, \mathbb{Z})$ with the negative of the intersection form as quadratic form. Therefore the difference $\overline{\Psi}_C^{X, H_1} - \overline{\Psi}_C^{X, H_2}$ for two points $H_1, H_2 \in \overline{\mathbb{H}}_X$ should be given in terms of a theta function $\Theta_{L,c,c}^{h_1, h_2}$ (where as usual the small letter denotes the Poincaré dual class in $L$. In the next section we will define such theta functions and show that, at least for $f$ and $g$ rational points with $\mathbf{Q}(f) = \mathbf{Q}(g) = 0$, $\Theta_{L,c,c}^{f,g}$ has the properties of usual theta functions, enabling us to prove our main structural results.

2. Conjecture 5.1 from [G] can be interpreted as saying that Theorems 2.4 and 2.8 hold more generally in the case $H^1(X, \mathbb{Q}) = 0$, if by $H^2(X, \mathbb{Z}) + C/2$ we mean the set of all expressions $(2\xi + C)/2$ (counted with repetitions), with $\xi$ running through $H^2(X, \mathbb{Z})$.

3. Theta functions for indefinite lattices

The classical theta series associated to a positive definite lattice $L$ with quadratic form $Q$ and associated bilinear form $x \cdot y$ (see notation 3.1 below) is the sum

$$\Theta_L(\tau, x) := \sum_{\xi \in L} q^{Q(\xi)} e^{2\pi i \xi \cdot x}, \qquad (\tau \in \mathfrak{H}, \ x \in L_\mathbb{C} = L \otimes \mathbb{C}).$$

These theta series have well-known transformation properties. In particular the "Nullwert" $\Theta_L(\tau, 0)$ is a modular form of weight $r/2$, where $r$ is the rank of $L$. In this chapter we give a generalization to the case when $L$ is allowed to be indefinite. In particular we will consider the case when the type of $Q$ is $(r-1, 1)$. The theta series that we define depends not only on $L$, but also on two vectors $f, g \in L_{\mathbb{Q}}$ with $Q(f) = Q(g) = 0$, $f \cdot g < 0$. It is defined in a certain open subset of $\mathfrak{H} \times L_{\mathbb{C}}$ by the formula

$$\Theta_L^{f,g}(\tau, x) := \Big(\sum_{\xi \cdot f \geq 0 > \xi \cdot g} - \sum_{\xi \cdot g \geq 0 > \xi \cdot f}\Big) q^{Q(\xi)} e^{2\pi i \xi \cdot x} \tag{3.0.1}$$

($\xi$ running through $L$). We will also sketch a generalization to the case of type $(n-s, s)$ with $n - s \geq s > 1$, but this is not needed for our applications to the Donaldson invariants.

The main properties of $\Theta_L^{f,g}$ are proved by using an alternative definition. The idea of the construction is simple. When $L$ is the standard hyperbolic lattice $H$ generated by vectors $f, g$ with $Q(f) = 0$, $Q(g) = 0$, $f \cdot g = -1$ then $\Theta_L^{f,g}$ will be the function $F(\tau; u, v)$ studied in [Z] (but in fact going back to Kronecker, see [We]). If $L$ is the direct sum $H \oplus L_0$, with $L_0$ positive definite, then $\Theta_L^{f,g}$ is just the product of this function with the usual theta series of $L_0$. The general case is reduced to this by considering the sublattice $L' = \langle f, g \rangle \oplus \langle f, g \rangle^\perp$ of $L$ and averaging over cosets of $L'$ in $L$.

**3.1. The function F.** The main building block for the construction of theta functions will be the function $F(\tau; u, v) : \mathfrak{H} \times \mathbb{C}^2 \longrightarrow \mathbb{C}$, studied in [Z], which is defined for $0 < -\Re(u)/\Im(\tau) < 2\pi$, $0 < \Re(v)/\Im(\tau) < 2\pi$ by the formula

$$F(\tau; u, v) := \sum_{n \geq 0, m > 0} q^{nm} e^{-nu - mv} - \sum_{n > 0, m \geq 0} q^{nm} e^{nu + mv}. \tag{3.0.2}$$

(This formula is not given explicitly in [Z], but is easily proved to be equivalent to 3. below.) This funcion has the following properties (see [Z]).

1. $F(\tau; u, v)$ has a meromorphic continuation to $\mathfrak{H} \times \mathbb{C}^2$ with simple poles for $u$ or $v$ in $2\pi i(\mathbb{Z}\tau + \mathbb{Z})$ and no other poles.
2. $F(\tau; u, v) = F(\tau; v, u) = -F(\tau; -u, -v)$.
3. For $|\Re(u)/\Im(\tau)| < 2\pi$, $|\Re(v)/\Im(\tau)| < 2\pi$ it has the Fourier expansion

$$F(\tau; u, v) := -\frac{1}{1 - e^u} + \frac{1}{1 - e^{-v}} - 2 \sum_{n > 0, m > 0} \sinh(nu + mv) q^{nm}. \tag{3.0.3}$$

4. $F\Big(\frac{a\tau + b}{c\tau + d}; \frac{u}{c\tau + d}, \frac{v}{c\tau + d}\Big) = (c\tau + d) e^{\frac{cuv/2\pi i}{c\tau + d}} F(\tau; u, v)$ for all $\begin{pmatrix} a & b \\ c & d \end{pmatrix} \in SL(2, \mathbb{Z})$.
5. $F(\tau; u + 2\pi i(n\tau + s), v + 2\pi i(m\tau + t)) = q^{-nm} e^{-nu - mv} F(\tau; u, v)$ for all $n, m, s, t \in \mathbb{Z}$.
6. $F(\tau; u, v) = \dfrac{\eta(\tau)^3 \theta_{11}(\tau, u + v)}{\theta_{11}(\tau, u) \theta_{11}(\tau, v)}$.
7. $F(\tau; u, v) = \dfrac{u + v}{uv} \exp\Big(\sum_{k > 0} \dfrac{2}{k!} [u^k + v^k - (u + v)^k] G_k(\tau)\Big)$.

**3.2. Definition of the theta functions.**

**Notation 3.1.** For us a *lattice* is a free $\mathbb{Z}$ module $L$ together with a quadratic form $Q : L \to \frac{1}{2}\mathbb{Z}$, such that the associated bilinear form $x \cdot y := Q(x + y) - Q(x) - Q(y)$ is nondegenerate and $\mathbb{Z}$-valued. The extensions of the quadratic and bilinear form to $L_{\mathbb{C}} := L \otimes \mathbb{C}$, or more generally $L_R := L \otimes_{\mathbb{Z}} R$ for any $\mathbb{Z}$-module $R$, are denoted in the same way. Let $L^\vee := \{\mu \in L_{\mathbb{Q}} \mid \mu \cdot L \subseteq \mathbb{Z}\}$. $L$ is unimodular if and only if $L = L^\vee$. The *type* of $L$ is the pair $(r-s, s)$, where $r$ is the rank of $L$ and $s$ the largest rank of a sublattice of $L$ on which $Q$ is negative definite, and the *signature* $\sigma(L)$ is the number $r - 2s$.

From now on we assume that $s = 1$, i.e., that $L$ has type $(r-1, 1)$. Then the set of vectors $f \in L_\mathbb{R}$ with $Q(f) < 0$ has two components, the scalar product of any two such vectors being negative if they belong to the same component and positive if they belong to opposite components. We fix a vector $f_0 \in L_\mathbb{R}$ with $Q(f_0) < 0$ and let

$$C_L := \{ f \in L_\mathbb{R} \mid Q(f) < 0,\ f \cdot f_0 < 0 \}$$

("positive light-cone") be the component containing $f_0$; we further set

$$S_L := \{ f \in L \mid f \text{ primitive},\ Q(f) = 0,\ f \cdot f_0 < 0 \}.$$

(The $(r-1)$-dimensional hyperbolic space $C_L/\mathbb{R}_+$ is the natural domain of definition of automorphic forms with respect to $O(L)$, and $S_L$ is a set of representatives for the corresponding set of cusps $\{f \in L_\mathbb{Q} \mid Q(f) < 0,\ f \cdot f_0 < 0\}/\mathbb{Q}_+$.) For $f \in S_L$ put

$$D(f) := \{ (\tau, x) \in \mathfrak{H} \times L_\mathbb{C} \mid 0 < \Im(f \cdot x) < \Im(\tau) \},$$

and for $f \in C_L$ put $D(f) := \mathfrak{H} \times L_\mathbb{C}$.

**Notation 3.2.** For $t \in \mathbb{R}$ we put $\mu(t) := \begin{cases} 1, & t \geq 0, \\ 0, & t < 0. \end{cases}$

**Definition 3.3.** Let $f, g \in C_L \cup S_L$. For $(\tau, x) \in D(f) \cap D(g)$, we define the *theta function* of $L$ with respect to $(f, g)$ by the formula (3.0.1). More generally we put for $c, b \in L$ and $(\tau, x) \in D(f) \cap D(g)$

$$\Theta^{f,g}_{L,c,b}(\tau, x) := \sum_{\xi \in L + c/2} \left( \mu(\xi \cdot f) - \mu(\xi \cdot g) \right) q^{Q(\xi)} e^{2\pi i \xi \cdot (x + b/2)}, \tag{3.3.1}$$

so that $\Theta^{f,g}_{L,0,0} = \Theta^{f,g}_L$. It is clear that $\Theta^{f,g}_{L,c,b}$ depends up to sign only on the class of $c$ and $b$ in $L/2L$. We will later want to show that in case $f, g \in S_L$ the function $\Theta^{f,g}_{L,c,b}(\tau, x)$ has nice analytical properties. To see that $\Theta^{f,g}_{L,c,b}(\tau, x)$ is well-defined we have to see that the sum (3.3.1) converges absolutely and locally uniformly on $D(f) \cap D(g)$.

*Case 1: $f, g \in S_L$.* We check the convergence only for (3.0.1) (for (3.3.1) this is analogous). Let $N := -f \cdot g$. It is enough to show the absolute convergence of (3.0.1) with $\xi$ running in $\frac{1}{N}(\langle f, g \rangle \oplus \langle f, g \rangle^\perp)$ instead of in $L$. So it is enough to check the absolute convergence of $\sum_{\eta \in \frac{1}{N}\langle f,g \rangle^\perp} q^{Q(\eta)} e^{2\pi i \eta \cdot x}$ (which holds even on $\mathfrak{H} \times L_\mathbb{C}$ as $\langle f, g \rangle^\perp$ is negative definite), and that of

$$\sum_{a=1}^{\infty} \sum_{b=0}^{\infty} q^{ab/N} e^{2\pi i (af - bg) \cdot x / N}, \qquad \sum_{b=1}^{\infty} \sum_{a=0}^{\infty} q^{ab/N} e^{2\pi i (-af + bg) \cdot x / N},$$

which converge absolutely on $D(f) \cap D(g)$. We shall see in 3.4 that $\Theta^{f,g}_L$ and $\Theta^{f,g}_{L,c,b}$ have meromorphic continuations to $\mathfrak{H} \times L_\mathbb{C}$.

*Case 2: $f, g \in C_L$.* Then $\Theta^{f,g}_{L,c,b}(\tau, x)$ depends only on the classes of $f$ and $g$ modulo $\mathbb{R}^+$. To check the absolute and locally uniform convergence of (3.3.1) on $\mathfrak{H} \times L_\mathbb{C}$, we write $\xi \in L_\mathbb{R}$ as $af + bg + \xi_\perp$ with $a, b \in \mathbb{R}$ and $\xi_\perp \in \langle f, g \rangle_\mathbb{R}^\perp$. Since then $Q(\xi) = Q(af + bg) + Q(\xi_\perp)$, it is enough to show that there is an $C > 0$ such that for $\xi = af + bg \in \langle f, g \rangle_\mathbb{R}$ with $(\xi \cdot f)(\xi \cdot g) \leq 0$ we have $Q(\xi) > C(a^2 + b^2)$. On $\langle f, g \rangle$ the quadratic form has type $(1, 1)$, therefore we have $(f \cdot g)^2 > 4Q(f)Q(g)$. The inequality $(af \cdot g + 2bQ(g))(2aQ(f) + bf \cdot g) \leq 0$ can be rewritten $\frac{(f \cdot g)^2 + 4Q(f)Q(g)}{f \cdot g} ab \geq -2Q(f) a^2 - 2Q(g) b^2$. Therefore we get

$$Q(\xi) = Q(f) a^2 + (f \cdot g) ab + Q(g) b^2 \geq \frac{(f \cdot g)^2 - 4Q(f)Q(g)}{(f \cdot g)^2 + 4Q(f)Q(g)} \left( -Q(f) a^2 - Q(g) b^2 \right),$$

and, using that $Q(f) < 0$, $Q(g) < 0$, the result follows. Note that the argument also shows that in case $f, g \in C_L$ the sum (3.3.1) makes sense as a formal power series: for each integer $k$ it contains only a finite number of summands $q^{Q(\xi)} e^{2\pi i \xi \cdot (x + b/2)}$ with $Q(\xi) \leq k$.

*Case 3*: $f \in S_L, g \in C_L$. To see that (3.3.1) converges, we note that $\Theta_{L,c,b}^{g,f}(\tau, x) = \Theta_{L,c,b}^{h,f}(\tau, x) + \Theta_{L,c,b}^{g,h}(\tau, x)$ if $h \in C_L$ and any two of the three theta series converge. By the absolute convergence of (3.3.1) for $f, g \in C_L$ we can therefore assume that $g \in L$. We split the sum (3.3.1) into two parts, the first consisting of the summands with $\xi \cdot f \neq 0$ and the second consisting of those with $\xi \cdot f = 0$. The first sum converges for $|\Im(f \cdot x)/\Im(\tau)| < 1$. The second sum can be rewritten for $(\tau, x) \in D(f)$ as

$$\sum_{n \geq 0} e^{2\pi i n f \cdot (x + b/2)} \sum_{\substack{\xi \cdot f = 0 \\ f \cdot g \leq \xi \cdot g < 0}} q^{Q(\xi)} e^{2\pi i \xi \cdot (x + b/2)}$$

$$= \frac{1}{1 - e^{2\pi i f \cdot (x + b/2)}} \sum_{t \in P_0} \sum_{\xi \in \langle f, g \rangle^\perp} q^{Q(\xi + t)} e^{2\pi i (\xi + t) \cdot (x + b/2)},$$

for a suitable finite set $P_0$. Thus the result follows as $\langle f, g \rangle^\perp$ is positive definite. This also shows that $\Theta_{L,c,b}^{f,g}(\tau, x)$ for $f \in S_L$ and $g \in C_L$ has a meromorphic extension to $|\Im(f \cdot x)/\Im(\tau)| < 1$, and is given there by the Fourier expansion

$$\Theta_{L,c,b}^{f,g}(\tau, x) := \sum_{\xi \cdot f \neq 0} \left( \mu(\xi \cdot f) - \mu(\xi \cdot g) \right) q^{Q(\xi)} e^{2\pi i \xi \cdot (x + b/2)} \qquad (3.3.2)$$

$$+ \frac{1}{1 - e^{2\pi i f \cdot (x + b/2)}} \sum_{\substack{\xi \cdot f = 0 \\ f \cdot g \leq \xi \cdot g < 0}} q^{Q(\xi)} e^{2\pi i \xi \cdot (x + b/2)},$$

the sums running through $\xi \in L + c/2$. The argument also shows that (3.3.2) makes sense as a formal power series, i.e. for every $k \in \mathbb{Z}_{>0}$ there are only finitely many $\xi$ in the sum with $Q(\xi) < k$.

**Remark 3.4.** For $f, g, h \in C_L \cup S_L$ and $(\tau, x) \in D(f) \cap D(g) \cap D(h)$ we have the cocycle condition $\Theta_{L,c,b}^{f,g}(\tau, x) + \Theta_{L,c,b}^{g,h}(\tau, x) = \Theta_{L,c,b}^{f,h}(\tau, x)$.

*Proof.* This is immediate from the definitions. □

Note that the set $D(f) \cap D(g) \cap D(h)$ is always nonempty. Therefore the cocycle condition continues to hold after meromorphic extension of the $\Theta_{L,c,b}^{f,g}(\tau, x)$.

3.3. **Jacobi forms.** We briefly recall the notion of a Jacobi form in the form in which we need it. (For more details in the one variable case see [E-Z]. In the one-variable case $Q(x)$ is just $mx^2$, where $m$ is the index.)

**Definition 3.5.** Let $L$ be a lattice of rank $r$, and denote the quadratic form on $L$ by $Q: L \longrightarrow \frac{1}{2}\mathbb{Z}$ as usual. Denote by $M_L$ the set of meromorphic maps $f: \mathfrak{H} \times L_\mathbb{C} \to \mathbb{C}$. For $v = (\lambda, \mu) \in L_{\mathbb{R}^2}$ and for $A = \begin{pmatrix} a & b \\ c & d \end{pmatrix} \in \Gamma$, and $k \in \mathbb{Z}$ we define maps $|v: M_L \longrightarrow M_L$ and $|_k A: M_L \longrightarrow M_L$ by putting

$$f|v(\tau, x) \quad := \quad q^{Q(\lambda)} \exp(2\pi i(\lambda \cdot (x + \mu/2))f(\tau, x + \lambda\tau + \mu) \quad (3.5.1)$$

$$f|_k A(\tau, x) \quad := \quad (c\tau + d)^{-k} \exp\left(-2\pi i \frac{cQ(x)}{c\tau + d}\right) f\left(\frac{a\tau + b}{c\tau + d}, \frac{x}{c\tau + d}\right). \quad (3.5.2)$$

We view elements $v = (v_1, v_2)$, $w = (w_1, w_2) \in L_{\mathbb{R}^2}$ as row vectors with entries in $L_\mathbb{R}$ and denote by $\langle v, w \rangle := v_1 \cdot w_2 - v_2 \cdot w_1$ the corresponding "determinant" and by $vA$ the application of $A \in SL(2, \mathbb{R})$ to $v$. Then $|v$ and $|_k A$ have the following compatibility properties

**Remark 3.6.**   1. $|_k A$ defines an action of the group $SL(2, \mathbb{R})$ on $M_L$.
   2. $(f|v)|w = e^{\pi i \langle v, w \rangle} f|(v + w) = e^{2\pi i \langle v, w \rangle}(f|w)|v$,
   3. $(f|v)|_k A = (f|_k A)|vA$.

*Proof.* The proof is elementary and is completely analogous to that of Theorem 1.4 in [E-Z]. □

In the situation of Definition 3.5, we call a function $f \in M_L$ a *holomorphic Jacobi form* of weight $k$ with respect to $(\Lambda, \Gamma)$ if

   1. $f|(\lambda, \mu) = (-1)^{\lambda \cdot \mu} f$ for all $(\lambda, \mu) \in L^2$.
   2. $f|_k A = f$ for all $A \in \Gamma$.
   3. $f$ is holomorphic in $\mathfrak{H} \times L_\mathbb{C}$ and holomorphic in the cusps (in a suitable sense).

A function $f \in M_L$ satisfying only 1. and 2. will be called a *meromorphic Jacobi form* of weight $k$. In this paper we use only meromorphic Jacobi forms and hence do not explain the holomorphy condition at the cusps in 3.

**Remark 3.7.** It is evident from the definitions that for a Jacobi form $f$ of weight $k$ for a lattice $L$ and $\Gamma \subset SL(2, \mathbb{Z})$ the function $(\tau, z) \mapsto f(\tau, yz)$ on $\mathfrak{H} \times \mathbb{C}$ will be for every $y \in L$ a Jacobi form of weight $k$ and index $Q(y)$ for $\Gamma$ in the sense of [E-Z].

**Example 3.8.** Let $L$ be a positive definite lattice of rank $r$, and let $T, V, S$ be as in (2.5.4). The theta function $\Theta_L : \mathfrak{H} \times L_\mathbb{C} \to \mathbb{C}$ of $L$ is given by

$$\Theta_L(\tau, x) := \sum_{\xi \in L} q^{Q(\xi)} e^{2\pi i \xi \cdot x} = \sum_{\xi \in L} 1|(\xi, 0).$$

We obviously have $\Theta_L|(\lambda, \mu) = (-1)^{\lambda \cdot \mu} \Theta_L$ for all $(\lambda, \mu) \in L \times L$, and $\theta^{-r}\Theta_L|V = \Theta_L$, $\theta^{-r}\Theta_L|T = (\theta_{01}^0)^{-r}\Theta_L|(0, w/2)$ for a characteristic vector $w$ of $L$. Note that

$$\Theta_L = \sum_{l \in P} \Theta_{L'}|(l, 0)$$

for any sublattice $L' \subset L$ and $P$ a system of representatives of $L/L'$.

Let $N$ be the index of $L$ in $L^\vee$, and let $P$ be a system of representatives for $L^\vee$ modulo $L$. Then it is a standard fact that
$$\left(\theta^{-r}\Theta_L\right)\big|S = \frac{1}{\sqrt{N}}\sum_{t\in P}(\theta^{-r}\Theta_L)\big|(t,0).$$

In particular, if $L$ is unimodular, then $\theta^{-r}\Theta_L$ is a (meromorphic) Jacobi form of weight 0 for $\Gamma_\theta$, and if $L$ is in addition even, it is a Jacobi form of weight 0 for $SL(2,\mathbb{Z})$.

### 3.4. Properties of $\Theta_L^{f,g}$ for $f, g \in S_L$.

**Theorem 3.9.** *Let $L$ be a unimodular lattice of type $(r-1,1)$; let $f, g \in S_L$ and $c, b \in L$. Then*

1. *$\Theta_L^{f,g}$ and $\Theta_{L,c,b}^{f,g}$ have meromorphic extensions to $\mathfrak{H} \times L_\mathbb{C}$.*
2. *$\Theta_{L,c,b}^{f,g}(\tau, x) = -(-1)^{c \cdot b}\Theta_{L,c,b}^{f,g}(\tau, -x)$.*
3. *For $|\Im(f \cdot x)/\Im(\tau)| < 1$, $|\Im(g \cdot x)/\Im(\tau)| < 1$ they have the Fourier developments*

$$\Theta_L^{f,g}(\tau, x) = \frac{1}{1 - e^{2\pi i f \cdot x}}\sum_{\substack{\xi \cdot f = 0 \\ f \cdot g \leq \xi \cdot g < 0}} q^{Q(\xi)}e^{2\pi i \xi \cdot x} - \frac{1}{1 - e^{2\pi i g \cdot x}}\sum_{\substack{\xi \cdot g = 0 \\ f \cdot g \leq \xi \cdot f < 0}} q^{Q(\xi)}e^{2\pi i \xi \cdot x}$$
$$+ 2\sum_{\xi \cdot f > 0 > \xi \cdot g} q^{Q(\xi)}\sinh(2\pi i \xi \cdot x), \tag{3.9.1}$$

$$\Theta_{L,c,b}^{f,g}(\tau, x) = \frac{1}{1 - e^{2\pi i f \cdot (x+b/2)}}\sum_{\substack{\xi \cdot f = 0 \\ f \cdot g \leq \xi \cdot g < 0}} q^{Q(\xi)}e^{2\pi i \xi \cdot (x+b/2)}$$
$$- \frac{1}{1 - e^{2\pi i g \cdot (x+b/2)}}\sum_{\substack{\xi \cdot g = 0 \\ f \cdot g \leq \xi \cdot f < 0}} q^{Q(\xi)}e^{2\pi i \xi \cdot (x+b/2)} \tag{3.9.2}$$
$$+ 2\sum_{\xi \cdot f > 0 > \xi \cdot g} q^{Q(\xi)}\sinh\left(2\pi i \xi \cdot (x + b/2)\right),$$

*where in (3.9.1) and (3.9.2) the $\xi$ run through $L$ and $L + c/2$ respectively.*

4. *$\Theta_L^{f,g}/\theta^{\sigma(L)}$ is a meromorphic Jacobi form of weight 1 for $(L, \Gamma_\theta)$.*
5. *For all $\lambda, \mu$ in $L$ and for any characteristic vector $w$ of $L$ we have*

$$\Theta_{L,c,b}^{f,g}|(\lambda,\mu) = (-1)^{c\mu - b\lambda + \lambda\mu}\Theta_{L,c,b}^{f,g},$$
$$(\Theta_{L,c,b}^{f,g}/\theta^{\sigma(L)})|_1 S = 1^{-b \cdot c/4}\Theta_{L,b,c}^{f,g}/\theta^{\sigma(L)},$$
$$(\Theta_{L,c,b}^{f,g}/\theta^{\sigma(L)})|_1 T = 1^{3Q(c)/4 - c \cdot w/4}\Theta_{L,c,b-c+w}^{f,g}/(\theta_{01}^0)^{\sigma(L)},$$
$$(\Theta_{L,c,b}^{f,g}/\theta^{\sigma(L)})|_1 V = 1^{Q(c)/2}\Theta_{L,c,b}^{f,g}/\theta^{\sigma(L)},$$
$$(\Theta_{L,c,b}^{f,g}/\theta^{\sigma(L)})|_1 W = 1^{-Q(c)/4 - c \cdot b/4}\Theta_{L,w-c+b,c}^{f,g}/(\theta_{10}^0)^{\sigma(L)}.$$

*Proof.* The main idea in the proof is to give an alternative definition of the functions $\Theta_L^{f,g}$ and $\Theta_{L,c,b}^{f,g}$, which allows us to relate them to the function $F(\tau;u,v)$ from Section 3.1 and theta functions of positive definite lattices.

In this proof $|A$ with $A \in SL(2,\mathbb{Z})$ always stands for $|_1 A$. If $H$ is the hyperbolic lattice of type $(1,1)$, and $f, g$ are generators of $H$ with $Q(f) = Q(g) = 0$, $f \cdot g = -1$, then we get immediately from the formulas (3.0.1) and (3.0.2) that

$$\Theta_H^{f,g}(\tau, x) = F(\tau; -2\pi i f \cdot x, 2\pi i g \cdot x)$$

for all $(\tau, x) \in \mathfrak{H} \times H_{\mathbb{C}}$ with $0 < \Im(f \cdot x) < \Im(\tau)$, $0 < \Im(g \cdot x) < \Im(\tau)$, and by the results cited in Section 3.1 this shows properties 1., 2., 3. and 4. for $\Theta_H^{f,g}(\tau, x)$.

Similarly let now $f$ and $g$ be vectors generating a lattice $L$ with quadratic form given by $Q(f) = Q(g) = 0$, $f \cdot g = -N \in \mathbb{Z}_{<0}$. Then we also see from the results of Section 3.1 that

$$\Theta_{f,g}(\tau, x) := F(N\tau; -2\pi i f \cdot x, 2\pi i g \cdot x)$$

fulfills $\Theta_{f,g}|(\lambda, \mu) = (-1)^{\lambda \cdot \mu} \Theta_{f,g}$ for all $(\lambda, \mu) \in L \times L^{\vee}$. For a system $R$ of representatives of $L$ modulo $L^{\vee}$ we also get immediately that

$$\Theta_{f,g}|S = \frac{1}{N} \sum_{t \in R} \Theta_{f,g}|(t, 0),$$

and $\Theta_{f,g}|T = \Theta_{f,g}$

Now let $L$ be arbitrary of type $(r-1, 1)$ and $f, g \in S_L$ with $f \cdot g = -N \in \mathbb{Z}_{<0}$. Let $L_0 := \langle f, g \rangle \oplus \langle f, g \rangle^{\perp}$. Let $P$ be a system of representatives for $L$ modulo $L_0$. For each $x \in L_{\mathbb{C}}$ we denote by $x_{f,g}$ and $x_{\perp}$ the orthogonal projections to $\langle f, g \rangle_{\mathbb{C}}$ and $\langle f, g \rangle_{\mathbb{C}}^{\perp}$. We write $\Theta_{\perp}$ for $\Theta_{\langle f,g \rangle^{\perp}}$. For $\tau \in \mathfrak{H}$, $x \in L_{\mathbb{C}}$, $t \in L_{\mathbb{Q}}$ we denote

$$(\Theta_{f,g} \boxtimes \Theta_{\perp})(\tau, x) := \Theta_{f,g}(\tau, x_{f,g}) \cdot \Theta_{\perp}(\tau, x_{\perp}),$$
$$((\Theta_{f,g} \boxtimes \Theta_{\perp})|(t, 0))(\tau, x) := (\Theta_{f,g}|(t_{f,g}, 0))(\tau, x_{f,g}) \cdot (\Theta_{\perp}|(t_{\perp}, 0))(\tau, x_{\perp}).$$

Note that for $t \in L$ the function $(\Theta_{f,g} \boxtimes \Theta_{\perp})|(t, 0)$ depends only class of $t$ in $L/L_0$. We define the function $\widetilde{\Theta}_L^{f,g}$ on $\mathfrak{H} \times L_{\mathbb{C}}$ by

$$\widetilde{\Theta}_L^{f,g} := \sum_{t \in P} (\Theta_{f,g} \boxtimes \Theta_{\perp})|(t, 0), \tag{3.9.3}$$

and for $c, b \in L$ we write $\widetilde{\Theta}_{L,c,b}^{f,g} := \widetilde{\Theta}_L^{f,g}|(c/2, 0)|(0, b/2)$.

*Claim.* For $0 < \Im(f \cdot x) < \Im(\tau)$, $0 < \Im(g \cdot x) < \Im(\tau)$ we have $\Theta_{L,c,b}^{f,g}(\tau, x) = \widetilde{\Theta}_{L,c,b}^{f,g}(\tau, x)$ (in particular we get the meromorphic continuation of $\Theta_{L,c,b}^{f,g}$ to $\mathfrak{H} \times L_{\mathbb{C}}$).

To prove the claim we can choose a system of representatives $P$ of $L/L_0$ such that all $t \in P$ satisfy $-f \cdot g > (t + c/2) \cdot f \geq 0$, $-f \cdot g > (t + c/2) \cdot g \geq 0$. Then we get for every $\xi \in L_0$ that $\xi \cdot f \geq 0$ if and only if $(\xi + t + c/2) \cdot f \geq 0$ and $\xi \cdot g \geq 0$ if and only if $(\xi + t + c/2) \cdot g \geq 0$. Using this the formula (3.5.1) gives for $t \in P$ and $(\tau, x) \in D(f) \cap D(g)$ that

$$((\Theta_{f,g} \boxtimes \Theta_{\perp})|(t, 0)|(c/2, 0)|(0, b/2))(\tau, x) = \sum_{\xi \in L_0 + c/2 + t} (\mu(\xi \cdot f) - \mu(\xi \cdot g)) q^{Q(\xi)} e^{2\pi i \xi \cdot (x + b/2)},$$

and the claim follows by summing over $t \in P$. We will in future also write $\Theta_L^{f,g}$ and $\Theta_{L,c,b}^{f,g}$ for their meromorphic extensions $\widetilde{\Theta}_L^{f,g}$ and $\widetilde{\Theta}_{L,c,b}^{f,g}$

2. By Section 3.1 we have $\Theta_{f,g}(\tau, -x_{f,g}) = -\Theta_{f,g}(\tau, x_{f,g})$ and by Example 3.8 $\Theta_\perp(\tau, -x_\perp) = \Theta_\perp(\tau, x_\perp)$. Thus the definition of $\widetilde{\Theta}_L^{f,g}$ gives that $\Theta_L^{f,g}(\tau, -x) = -\Theta_L^{f,g}(\tau, x)$. Therefore

$$\Theta_{L,c,b}^{f,g}(\tau, -x) = \left(\Theta_L^{f,g}|(c/2, 0)|(0, b/2)\right)(\tau, -x)$$
$$= -\left(\Theta_L^{f,g}|(-c/2, 0)|(0, -b/2)\right)(\tau, x) = -(-1)^{b \cdot c}\Theta_{L,c,b}^{f,g}(\tau, x).$$

3. We can assume $b = 0$. By definition we have on $D(f) \cap D(g)$:

$$\Theta_{L,c,0}^{f,g}(\tau, x) = \left(\sum_{\substack{\xi \cdot f = 0 \\ \xi \cdot g < 0}} - \sum_{\substack{\xi \cdot g = 0 \\ \xi \cdot f < 0}}\right) q^{Q(\xi)} e^{2\pi i \xi \cdot x} + 2 \sum_{\xi \cdot f > 0 > \xi \cdot g} q^{Q(\xi)} \sinh(2\pi i \xi \cdot x),$$

the sums running over $\xi \in L + c/2$. The last sum converges obviously for $|\Im(f \cdot x)/\Im(\tau)| < 1$, $|\Im(g \cdot x)/\Im(\tau)| < 1$, and on $D(f) \cap D(g)$ we have

$$\sum_{\substack{\xi \cdot f = 0 \\ \xi \cdot g < 0}} q^{Q(\xi)} e^{2\pi i \xi \cdot x} = \frac{1}{1 - e^{2\pi i f \cdot x}} \sum_{\substack{\xi \cdot f = 0 \\ f \cdot g \leq \xi \cdot g < 0}} q^{Q(\xi)} e^{2\pi i \xi \cdot x}.$$

The sum on the right hand side converges on $\mathfrak{H} \times L_{\mathbb{C}}$, as it can be rewritten as

$$\sum_{\substack{t \in P + c/2 \\ t \cdot f = 0}} \left(\sum_{\xi \in \langle f, g \rangle^\perp} q^{Q(\xi)} e^{2\pi i \xi \cdot x}\right) \bigg|(t, 0)$$

for a suitable system $P$ of representatives of $L/L_0$, and the inner sum converges absolutely because $\langle f, g \rangle_\perp$ is positive definite.

To prove the theorem it is now enough to show 5. Let $(\lambda, \mu) \in L^2$. Then

$$\Theta_{L,c,b}^{f,g}|(\lambda, \mu) = \sum_{t \in P}(\Theta_{f,g} \boxtimes \Theta_\perp)|(t, 0)|(c/2, 0)|(0, b/2)|(\lambda, \mu)$$
$$= \sum_{t \in P}(-1)^{c \cdot \mu - b \cdot \lambda + \lambda \cdot \mu}(\Theta_{f,g} \boxtimes \Theta_\perp)|(0, \mu)|(\lambda + t, 0)|(c/2, 0)|(0, b/2)$$
$$= \sum_{t' \in P'}(-1)^{c \cdot \mu - b \cdot \lambda + \lambda \cdot \mu}(\Theta_{f,g} \boxtimes \Theta_\perp)|(t', 0)|(c/2, 0)|(0, b/2)$$
$$= (-1)^{c \cdot \mu - b \cdot \lambda + \lambda \cdot \mu}\Theta_{L,c,b}^{f,g},$$

where $P'$ is another system of representatives. Let $w$ be a characteristic vector of $L$. Then

$$(\Theta_{L,c,b}^{f,g}/\theta^{\sigma(L)})|T = \sum_{t \in P}\left((\Theta_{f,g} \boxtimes \Theta_\perp)/\theta^{\sigma(L)}\right)|(t, 0)|(c/2, 0)|(0, b/2)|T$$
$$= \sum_{t \in P}\left((\Theta_{f,g} \boxtimes \Theta_\perp)/\theta^{\sigma(L)}\right)|T|(t, t)|(c/2, c/2)|(0, b/2)$$
$$= (\theta_{01}^0)^{-\sigma(L)}\sum_{t \in P}(\Theta_{f,g} \boxtimes \Theta_\perp)|(0, w/2)|(t, t)|(c/2, c/2)|(0, b/2)$$
$$= (\theta_{01}^0)^{-\sigma(L)}\sum_{t \in P}(-1)^{t \cdot (t-w)}(\Theta_{f,g} \boxtimes \Theta_\perp)|(t, 0)|(0, w/2)|(c/2, c/2)|(0, b/2)$$
$$= 1^{3Q(c)/4 - c \cdot w/4}(\theta_{01}^0)^{-\sigma(L)}\Theta_{L,c,b-c+w}^{f,g}.$$

Applying $|T$ twice we get

$$(\Theta_{L,c,b}^{f,g}/\theta^{\sigma(L)})|V = 1^{Q(c)/2}(\Theta_{L,c,b}^{f,g}/\theta^{\sigma(L)}).$$

Let $R$ be a system of representatives for $L_0^\vee/L_0$ and let $N$ be the index of $L_0$ in $L$. Then

$$(\Theta_{L,c,b}^{f,g}/\theta^{\sigma(L)})|S = \sum_{t \in P}(\Theta_{f,g}|S \boxtimes (\Theta_\perp/\theta^{\sigma(L)})|S)|(0,-t)|(0,-c/2)|(b/2,0)$$

$$= \sum_{t \in P}\sum_{r \in R} \frac{1}{N}\theta^{-\sigma(L)}(\Theta_{f,g} \boxtimes \Theta_\perp)|(r,0)|(0,-t)|(0,-c/2)|(b/2,0)$$

$$= 1^{-b\cdot c/4}\sum_{r \in R}\left(\sum_{t \in P}\frac{1}{N}1^{-ir\cdot t}\theta^{-\sigma(L)}(\Theta_{f,g}\boxtimes\Theta_\perp)|(r,0)\right)|(b/2,0)|(0,c/2).$$

The inner sum is zero unless $r \in L^\vee = L$. Thus we get $(\Theta_{L,c,b}^{f,g}/\theta^{\sigma(L)})|S = 1^{-b\cdot c/4}\Theta_{L,b,c}^{f,g}/\theta^{\sigma(L)}$.

Finally we get

$$(\Theta_{L,c,b}^{f,g}/\theta^{\sigma(L)})|W = (\Theta_{L,c,b}^{f,g}/\theta^{\sigma(L)})|V^{-1}|T|S$$

$$= 1^{Q(c)/4-c\cdot w/4}\left((\theta_{01}^0)^{-\sigma(L)}\Theta_{L,c,b-c+w}^{f,g}\right)|S$$

$$= 1^{-Q(c)/4-b\cdot c/4}(\theta_{10}^0)^{-\sigma(L)}\Theta_{L,w-c+b,c}^{f,g}.$$

$\square$

**Remark 3.10.** We mention a (very) partial generalization of the theta functions to unimodular lattices $L$ of type $(r-s,s)$ with $r-s \geq s$. Let $F := (f_1,\ldots,f_s)$, $G := (g_1,\ldots,g_s)$, where the $f_j$, $g_j$ are primitive elements of $L$ with $Q(f_j) = Q(g_j) = 0$, $f_j \cdot g_j < 0$ and $f_j \cdot f_l = f_j \cdot g_l = g_j \cdot g_l = 0$ if $j \neq l$. On the set of $(\tau,x)$ with $0 < \Im(f_j \cdot x) < \Im(\tau)$, $0 < \Im(g_j \cdot x) < \Im(\tau)$ for all $j$ we define $\Theta_L^{F,G}$ and $\Theta_{L,c,b}^{F,G}$ for $c,b \in L$ by

$$\Theta_L^{F,G}(\tau,x) = \sum_{\xi \in L}\Big(\prod_{i=1}^s \big(\mu(f_i\cdot\xi) - \mu(g_i\cdot\xi)\big)\Big)q^{Q(\xi)}e^{2\pi i\xi\cdot x}, \tag{3.10.1}$$

$$\Theta_{L,c,b}^{F,G}(\tau,x) = \sum_{\xi \in L+c/2}\Big(\prod_{i=1}^s \big(\mu(f_i\cdot\xi) - \mu(g_i\cdot\xi)\big)\Big)q^{Q(\xi)}e^{2\pi i\xi\cdot(x+b/2)}. \tag{3.10.2}$$

Now let $f := f_s$, $g := g_s$, $F_s := (f_1,\ldots,f_{s-1})$, $G_s := (g_1,\ldots,g_{s-1})$, and let $L_0 = \langle f,g\rangle \oplus \langle f,g\rangle^\perp$ and $\Theta_\perp := \Theta_{\langle f,g\rangle^\perp}^{F_s,G_s}$. Then the same arguments as above show that $\Theta_L^{F,G}$ is defined inductively by the right hand side of (3.9.3) and that $\Theta_{L,c,b}^{F,G} = \Theta_L^{F,G}|(c/2,0)|(0,b/2)$. Furthermore 1., 2., 3. and 5. of Theorem 3.9 still hold, if we replace everywhere $\Theta_L^{f,g}$ by $\Theta_L^{F,G}$ and $\Theta_{L,c,b}^{f,g}$ by $\Theta_{L,c,b}^{F,G}$ and interpret 3. as an inductive definition of a Fourier development of $\Theta_L^{F,G}$, $\Theta_{L,c,b}^{F,G}$ on the set of $(\tau,x)$ with $|\Im(f_j\cdot x)/\Im(\tau)| < 1$, $|\Im(g_j\cdot x)/\Im(\tau)| < 1$ for all $j$, using the notations just introduced. We also have to replace $|_1$ by $|_s$, and instead of 4. we obtain that $\Theta_L^{F,G}/\theta^{\sigma(L)}$ is a meromorphic Jacobi form of weight $s$ for $L$ and $\Gamma_\theta$. The proofs are staightforward generalizations of the proof of Theorem 3.9.

### 3.5. The structure theorem for the theta functions.

**Notation 3.11.** During this section let $L$ be a unimodular lattice of type $(r-1,1)$, let $c \in L$, and let $f,g \in S_L$. Let $x \in L$, and assume $f \cdot x \neq 0$ if $c \cdot f$ is even and $g \cdot x \neq 0$ if $c \cdot g$ is even. Then we put for $(\tau, z) \in \mathfrak{H} \times \mathbb{C}$

$$\varphi_{L,c}^{f,g}(\tau, x \cdot z) := 1^{-3/8c \cdot c} \frac{2\theta(\tau)^{-\sigma(L)}}{f(\tau)} \Theta_{L,c,c}^{f,g}\left(\tau, \frac{xz}{2\pi i f(\tau)}\right) e^{2Q(x)G(\tau)z^2/f(\tau)^2}.$$

When in future we write $\varphi_{L,c}^{f,g}(\tau, x \cdot z)$ we will assume implicitly that the above conditions on $f, g, c, x$ are fulfilled.

It is well known that for a (one-variable) Jacobi form $\phi(\tau, z)$ of weight $k$ and index $m$ for $\Gamma \subset SL(2, \mathbb{Z})$ the coefficient $w_n(\tau)$ of $z^n$ in the Taylor development

$$\phi(\tau, z) e^{-8\pi^2 m z^2 G_2(\tau)} := \sum_{n \geq 0} w_n(\tau) z^n$$

is a modular form of weight $k + n$ for $\Gamma$. This follows readily from the transformation behaviour

$$G_2\left(\frac{a\tau + b}{c\tau + d}\right) = (c\tau + d)^2 G_2(\tau) - \frac{c(c\tau + d)}{4\pi i}, \qquad \begin{pmatrix} a & b \\ c & d \end{pmatrix} \in SL(2, \mathbb{Z}) \tag{3.11.1}$$

of the quasi-modular form $G_2$. We show a similar result for $\varphi_{L,c}^{f,g}(\tau, x \cdot z)$, giving however a much more precise description of the coefficients.

**Definition 3.12.** Let $m = m(L, f, g) := \min\{w \cdot w \mid w \in L \text{ characteristic}, w \cdot f \geq 0 \geq w \cdot g\}$. Note that $(\sigma(L) - m)/8$ is an integer. Let

$$B(L, f, g) := \{w \in L \text{ characteristic} \mid w \cdot w < \sigma(L), (w \cdot f)(w \cdot g) \leq 0\},$$
$$B_i(L, f, g) := \{w \in B(L, f, g) \mid w \cdot f < 0 < w \cdot g\},$$
$$B_f(L, f, g) := \{w \in B(L, f, g) \mid 0 = w \cdot f > w \cdot g \geq 2f \cdot g\},$$
$$B_g(L, f, g) := \{w \in B(L, f, g) \mid 0 = w \cdot g > w \cdot f \geq 2f \cdot g\}.$$

The elements $w \in B(L, f, g)$ are called *basic classes* for $(L, f, g)$. Usually we will drop $(L, f, g)$ in the notation. Note that the sets $B_i$, $B_f$ and $B_g$ are finite. For $(\tau, x) \in \mathfrak{H} \otimes L_{\mathbb{C}}$ we put

$$\Omega_{L,c}^{f,g}(\tau, x) := -\sum_{w \in B_f} \frac{(-1)^{c \cdot (w+c)/2} q^{Q(w)/4} e^{-w \cdot x}}{1 - (-1)^{c \cdot f} e^{-2f \cdot x}} + \sum_{w \in B_g} \frac{(-1)^{c \cdot (w+c)/2} q^{Q(w)/4} e^{-w \cdot x}}{1 - (-1)^{c \cdot f} e^{-2g \cdot x}}$$
$$+ \sum_{w \in B_i} q^{Q(w)/4} \left((-1)^{c \cdot (w+c)/2} e^{-w \cdot x} - (-1)^{c \cdot (-w+c)/2} e^{w \cdot x}\right). \tag{3.12.1}$$

For $x \in L_{\mathbb{C}}$ we also write

$$O_{L,c}^{f,g}(x) := -\sum_{\substack{w \in B_f \\ w \cdot w = m}} \frac{(-1)^{c \cdot (w+c)/2} e^{-w \cdot x}}{1 - (-1)^{c \cdot f} e^{-2f \cdot x}} + \sum_{\substack{w \in B_g \\ w \cdot w = m}} \frac{(-1)^{c \cdot (w+c)/2} e^{w \cdot x}}{1 - (-1)^{c \cdot g} e^{-2g \cdot x}}$$
$$+ \sum_{\substack{w \in B_i \\ w \cdot w = m}} \left((-1)^{c \cdot (w+c)/2} e^{-w \cdot x} - (-1)^{c \cdot (-w+c)/2} e^{w \cdot x}\right), \tag{3.12.2}$$

so that $O_{L,c}^{f,g}(x)$ is the leading term of $\Omega_{L,c}^{f,g}(\tau, x)$ (coefficient of $q^{m/2}$) as $q \to 0$.

**Theorem 3.13.** 1. $\varphi_{L,c}^{f,g}(\tau, x \cdot z)$ has a Laurent development

$$\varphi_{L,c}^{f,g}(\tau, x \cdot z) = \sum_{n \geq -1} w_n(\tau) z^n,$$

where each $w_n(\tau)$ is a modular function for $\Gamma_u$. More precisely:

2.
$$w_n(\tau) = P_n\left(\frac{1}{U(\tau) - 2}\right) - 1^{(c \cdot c - 1 - n)/4} P_n\left(\frac{1}{-U(\tau) - 2}\right) + R_n(U(\tau)),$$

where $P_n(t)$ is a polynomial of degree $(\sigma(L) - m)/8$ in $t$ without constant term and $R_n(t)$ is a polynomial of degree $\leq (n+1)/2$.

3. The polynomials $P_n(t)$ are determined by the fact that $\sum_{n=-1}^{\infty} P_n\left(\frac{1}{\widetilde{U}(\tau) - 2}\right) z^n$ is the principal part in $(\widetilde{U}(\tau) - 2)$ of the development of

$$\frac{4\eta(2\tau)^2}{\theta_{10}^0(\tau)^{\sigma(L)} \eta(\tau)^4} \Omega_{L,c}^{f,g}\left(\tau, xz\frac{\eta(2\tau)^2}{\eta(\tau)^4}\right) \exp\left(Q(x)\left(8G_2(\tau) + 4e_1(\tau)\right)\frac{\eta(2\tau)^4}{\eta(\tau)^8} z^2\right)$$

as a Laurent series $\sum_{n=-1}^{\infty} \left(\sum_{k=-(\sigma(L)-m)/8}^{\infty} a_{n,k} (\widetilde{U}(\tau) - 2)^m\right) z^n$ in $z$ and $\widetilde{U}(\tau) - 2$.

4. The leading coefficients $a_n$ of $t^{(\sigma(L)-m)/8}$ in $P_n(t)$ are given by

$$\sum_{n=-1}^{\infty} a_n z^n = 2^{2-(\sigma(L)+3m)/4} O_{L,c}^{f,g}(xz) e^{-Q(x) z^2}.$$

*Proof.* 1. The properties of $f$ from 2.2 and Theorem 3.9 show that for $l \in \mathbb{Z}$ congruent to $2Q(c)$ modulo 2 the function $\theta(\tau)^{-\sigma(L)} f(\tau)^l \Theta_{L,c,c}^{f,g}(\tau, x)$ is a (meromorphic) Jacobi form of weight $l + 1$ for $\Gamma_u$ and $L$. Now let $x$ satisfy the conditions of notation 3.11. The Fourier development of $\Theta_{L,c,c}^{f,g}(\tau, x)$ from part 2. of Theorem 3.9, the fact that $e_3$ is a modular form for $\Gamma_\theta$ and (3.11.1) and (3.5.2) show that we get a Laurent development

$$\theta(\tau)^{-\sigma(L)} f(\tau)^l \Theta_{L,c,c}^{f,g}(\tau, xz) e^{-8\pi^2 Q(x) G(\tau) z^2} = \sum_{n \geq -1} v_n(\tau) z^n,$$

where each $v_n(\tau)$ has the transformation behaviour of a modular form of weight $l + n + 1$ for $\Gamma_u$, which has only poles at the zeroes of $\theta(\tau)^{\sigma(L)} f(\tau)^{-l}$. We see that $w_n(\tau) = 1^{-\frac{3}{8} c \cdot c} \frac{2 v_n(\tau)}{(2\pi i)^n f(\tau)^{1+l+n}}$, and the result follows as $f(\tau)^2$ is a modular form of weight 2 for $\Gamma_u$ ($v_n$ is zero unless $1 + l + n \in 2\mathbb{Z}$ by part 3. of Theorem 3.9).

2. As $U(\tau)$ defines an isomorphism $\mathfrak{H}/\Gamma_u \cup \{-1\} \cup \{1\} \cup \{\infty\} \to \mathbb{P}_1$ we see that $w_n(\tau)$ is a rational function in $U(\tau)$. The functions $\theta(\tau)$ and $f(\tau)$ are holomophic and nonzero on $\mathfrak{H}$, therefore $w_n(\tau)$ can only have poles in the cusps. As $U(\tau)$ sends the cusps $\infty$, $1$, $-1$ to $\infty$, $-2$, $2$ respectively, this shows that $w_n(\tau) = P_n\left(\frac{1}{U(\tau)-2}\right) + Q_n\left(\frac{1}{U(\tau)+2}\right) + R_n(U(\tau))$ for suitable polynomials $P_n$, $Q_n$ and $R_n$. As $\theta(\tau)$ is holomorphic and nonzero at $q = 0$, and the $q$-developments of $f(\tau)$ and $U(\tau)$ start in degree $1/8$ and $-1/4$, the degree of $R_n$ is at most $(n+1)/2$. To determine $P_n$ and its degree we apply $W \in SL(2, \mathbb{Z})$, which sends $\infty$ to $-1$. As $U|W = \widetilde{U}$, and the $q$-development of $\widetilde{U}(\tau) - 2$ starts in degree 1, the degree of $P_n$ is the order of pole of the $q$-development of $w_n|W(\tau)$. By definition the $w_n|W(\tau)$ are the coefficients in the Laurent development in $z$ of $\varphi_{L,c}^{f,g}(W\tau, x \cdot z)$, and by Theorem 3.9, and the results of Section 2.2 we get that

$$\varphi_{L,c}^{f,g}(W\tau, x \cdot z) := -\left(1^{c \cdot c/4} \frac{4\eta(2\tau)^2}{\theta_{10}^0(\tau)^{\sigma(L)}\eta(\tau)^4} \Theta_{L,w,c}^{f,g}\left(\tau, -xz\frac{\eta(2\tau)^2}{\pi i \eta(\tau)^4}\right)\right. \tag{3.13.1}$$

$$\left. \cdot \exp\left(Q(x)(8G_2(\tau) + 4e_1(\tau))\frac{\eta(2\tau)^4}{\eta(\tau)^8}\right)\right).$$

By 3.9 the lowest power of $q$ occuring in the Fourier development of $\Theta_{L,w,c}^{f,g}(\tau, x)$ is $q^{m/8}$. On the other hand the Fourier development of $(\theta_{10}^0)^{-\sigma(L)}$ starts with $q^{-\sigma(L)/8}$, and $G_2(\tau)$, $e_1(\tau)$, $\eta(2\tau)^2/\eta(\tau)^4$ are holomorphic and nonzero at $q = 0$. Therefore the degree of $P_n$ is $(\sigma(L) - m)/8$. The formula $Q_n(t) = 1^{(c \cdot c - n - 1)/4} P_n(-t)$ follows from the fact that $VW$ transports $\infty$ to $1$ and $f|V(\tau) = if(\tau)$, $U|V(\tau) = -U(\tau)$, $\Theta_{L,c,c}^{f,g}|V = 1^{c \cdot c/4}\Theta_{L,c,c}^{f,g}$, and therefore $\varphi_{L,c}^{f,g}(VW\tau, x \cdot z) = 1^{(c \cdot c - 1)/4}\varphi_{L,c}^{f,g}(W\tau, -ix \cdot z)$, and finally $U|V(\tau) - 2 = -(U(\tau) + 2)$.

The formula (3.13.1) shows that 3. holds if we replace $\Omega_{L,c}^{f,g}(\tau, x)$ by $-\left(1^{c \cdot c/4}\Theta_{L,w,c}^{f,g}(\tau, -\frac{x}{\pi i})\right)$. But by definition their $q$-developments are congruent modulo the ideal generated by $q^{\sigma(L)/8}$. Therefore the result follows.

To show 4. we just have to compute the leading terms of the expressions occuring in 3. So we use the congruences

$$\eta(\tau)^4/\eta(2\tau)^2 \equiv 1, \quad G_2(\tau) + e_1(\tau)/2 \equiv -1/8, \quad q^{\sigma(L)/8}(\theta_{10}^0)^{-\sigma(L)} \equiv 2^{-\sigma(L)}$$

modulo the ideal generated by $q$ and $\widetilde{U}(\tau) - 2 \equiv 64q$ modulo the ideal generated by $q^2$, and the result follows. $\square$

## 4. Application to Donaldson invariants

We want to apply the results about theta functions from Section 3 to get structural results for the Donaldson invariants of a simply connected 4-manifold with $b_+ = 1$. For this our lattice $L$ will be the lattice $H_2(X, \mathbb{Z})$ with the negative of the intersection form, i.e. $\mathbf{Q} = -2Q$, $\sigma(X) = -\sigma(L)$ and for $a, b$ in $H_2(X, \mathbb{Z})$ we have $a \cdot b = -AB$ where $AB$ is the intersection product of the Poincaré duals.

### 4.1. Extension of the Donaldson invariants.
The Donaldson invariant $\Phi_C^{X,H}$ is defined for $H \in \mathbb{H}_X$ with $H\xi \neq 0$ for all $\xi \in H^2(X, \mathbb{Z}) + C/2$. In (2.5.3) we used it to define new formal power series $\Psi_C^{X,H}$ and $\overline{\Psi}_C^{X,H}$. In this section we extend the definition of these invariants to arbitrary $H \in \overline{\mathbb{H}}_X$. To do this we apply the wall-crossing formulas Theorem 2.8 (we also use the notations from there) and use the theta functions of Section 3.

**Definition 4.1.** Fix $H \in \mathbb{H}_X$ with $H\xi \neq 0$ for all $\xi \in H^2(X, \mathbb{Z}) + C/2$. Let $M \in \overline{\mathbb{H}}_X$ and $x \in H^2(X, \mathbb{C})$, if $M \in \mathbb{S}_X$ and a primitive representative of $M$ has even intersection with $C$, we also assume that $xF \neq 0$. Denote again by $h, f, c \in L$ the Poincaré duals of $H, F, C$. Then we put

$$\Psi_C^{X,M}(x \cdot z, p^r) := \Psi_{C,r}^{X,H}(x \cdot z, p^r)$$
$$+ 1^{-3Q(c)/4} \mathrm{Coeff}_{u(\tau)^{r+1}} \left[ \frac{2\theta(\tau)^{\sigma(X)}}{f(\tau)} \Theta_{L,c,c}^{m,h}(\tau, \tfrac{xz}{2\pi i f(\tau)}) e^{2Q(x)G(\tau)z^2/f(\tau)^2} \right],$$
(4.1.1)
$$\overline{\Psi}_C^{X,M}(x \cdot z, t) := \sum_{r \geq 0} \Psi_C^{X,M}(x \cdot z, p^r) t^{r+1}.$$

Here we view the expression in square brackets as a formal Laurent series in $q^{1/8}$ and $z$ (i.e. it is a Laurent series in $z$, and the coefficient of $z^n$ is for every $n$ a Laurent series in $q^{1/8}$) using the formulas (3.3.1) and (3.3.2) (in case $M \in \mathbb{S}_X$) for $\Theta_{c,c}^{m,h}(\tau, x)$, which, as we have seen, also make sense as formal power series. We also define $\Phi_C^{X,M}(\alpha)$ for $\alpha \in A_*(X)$ by putting $\Phi_C^{X,M}(x^s p^r)$ to be $s!$ times the coefficient of $z^s$ in $\Psi_C^{X,M}(x \cdot z, p^r)$ and extending linearly (this is compatible with our previous definition of $\Phi_C^{X,M}$).

We check that $\overline{\Psi}_C^{X,M}(x \cdot z, t)$ is well defined. If $M \in \mathbb{H}_X$ with $M\xi \neq 0$ for all $\xi \in H^2(X, \mathbb{Z}) + C/2$, we have to see that the definition coincides with (2.5.1). By Theorem 2.8 we get (with the old definition) that

$$\Psi_C^{X,M}(x \cdot z, p^r) - \Psi_C^{X,H}(x \cdot z, p^r) = \sum_{M \cdot \xi > 0 > H \cdot \xi} 1^{C^2/8}(-1)^{(\xi - C/2)C} \mathrm{Coeff}_{u(\tau)^{r+1}} \left[ \Delta_\xi^X(\tau, x \cdot z) \right],$$

the sum running through $\xi \in L + c/2$. As $f(\tau)$ and $u(\tau)$ are power series in $q$ multiplied with $q^{1/8}$ and $q^{1/4}$ respectively we can replace $[\Delta_\xi^X(\tau, x \cdot z)]$ by $[\Delta_\xi^X(\tau, x \cdot z) - (-1)^{C^2} \Delta_{-\xi}^X(\tau, x \cdot z)]/2$, i.e. we get (4.1.1).

Now it follows immediately from the cocycle condition 3.4 that the definition of $\Psi_C^{X,M}(x \cdot z, p^r)$ above is independent of the choice of $H$.

**Remark 4.2.** The definition is motivated as follows:

1. For $M \in \mathbb{H}_X$ lying on a wall defined by a class of type $(C, d)$ the definition gives $\Phi_{C,d}^{X,M} = \Phi_{C,d}^{X,H} + \sum \delta_{\xi,d}^X / 2$ for $H$ in any chamber containing $L$ in its closure, and $\xi$ running through the classes of type $(C, d)$ through $M$ with $\xi H < 0$; in other words, we take the average over all chambers which contain $M$ in their closure.

2. If $F \in \mathbb{S}_X$, and its primitive representative has odd intersection with $C$, then $F$ lies in the closure of a unique chamber of type $(C, d)$. Our definition gives $\Phi_{C,d}^{X,F} := \Phi_{C,d}^{X,H}$ for $H$ in this chamber.

3. If $F \in \mathbb{S}_X$, and its primitive representative in $L$ has even intersection with $C$, then $F$ will in general lie on infinitely many walls defined by classes of type $(C, d)$, for every $d$. We would formally get from the definition (3.3.1) that

$$\Psi_C^{X,F}(x \cdot z, p^r) - \Psi_C^{X,H}(x \cdot z, p^r) =$$
$$\mathrm{Coeff}_{u(\tau)^{r+1}} \left[ \sum_{\xi \in H^2(X, \mathbb{Z}) + C/2} (\mu(F \cdot \xi) - \mu(H \cdot \xi)) 1^{C^2/8} (-1)^{(\xi - C/2)C} \Delta_\xi^X(\tau, x \cdot z) \right],$$

if the sum in square brackets converged as a formal power series in $z$ and $q$. Instead, we first make an analytic continuation. Therefore we can view the part of degree $d - 2r$ in $z$ of $\Psi_C^{X,M}(x \cdot z, p^r)$ as a "renormalized average" over the infinitely many chambers of type $(C, d)$ having $F$ in their closure.

Note that by (3.3.2) the function $\Psi_C^{X,M}(x \cdot z, p^r)$ will, in the case that both $CF$ and $C^2$ are even, usually be meromorphic in $z$ (with a simple pole at $z = 0$). This is the main reason why we introduced the notation $\Psi_C^{X,M}(x \cdot z, p^r)$ instead of $\Phi_C^{X,M}(e^{xz}p^r)$ for the Donaldson invariants.

With this definition the connection between the difference of the Donaldson invariants at period points $F, G \in \mathbb{S}_X$ on the one hand and the theta functions $\Theta_{L,c,c}^{f,g}$ and also the associated function $\varphi_{L,c}^{f,g}$ on the other becomes evident:

**Corollary 4.3.** *Let $F, G \in \mathbb{S}_X$ and $f$, $g$ Poincaré duals of representatives in $H^2(X, \mathbb{Z})$. Then*

$$\Psi_C^{X,F}(x \cdot z, p^r) - \Psi_C^{X,G}(x \cdot z, p^r) = \mathrm{Coeff}_{u(\tau)^{r+1}} \left[ \varphi_{L,c}^{f,g}(\tau, x \cdot z) \right].$$

*Proof.* This is straightforward from definition 4.1, 3.11 and the cocycle condition 3.4. □

4.2. **Blowup formulas.** The blowup formulas relate the Donaldson invariants of a 4-manifold $Y$ and $\widehat{Y} = Y \# \overline{\mathbb{P}}_2$. We have already used a small part of them in (2.5.2). In the case $b_+(Y) > 1$, when the invariants do not depend on the chamber structure, they have been shown in the most general form in [F-S1]. In [T1] they are shown also to hold in the case $b_+(Y) = 1$, if one takes the chamber structure into account (see also [K-L]). We cite only a weakened form, also avoiding the concept of related chambers.

**Theorem 4.4.** ([F-S1],[T1]) *There exist universal polynomials $B_k(t)$, $S_k(t) \in \mathbb{Q}[t]$ ($k = 0, 1, \ldots$) such that the following holds. Let $X$ be a simply connected 4-manifold, let $C \in H^2(X, \mathbb{Z})$ be not divisible by 2, and, in case $b_+(X) = 1$, let $M \in H^2(X, \mathbb{R})^+$ with $M\xi \neq 0$ for all $\xi \in H^2(X, \mathbb{Z}) + C/2$. Then for all $\alpha \in A_*(X)$ we have*

$$\Phi_C^{\widehat{X},M}(\alpha e^k) = \Phi_C^{X,M}(\alpha B_k(p)), \quad \Phi_{C+E}^{\widehat{X},M}(\alpha e^k) = \Phi_C^{X,M}(\alpha S_k(p)) \tag{4.4.1}$$

Here, as usual, $p \in A_*(X)$ denotes the class of a point. Using the notation of (2.5.3) and replacing $\alpha$ by $e^{xz}p^r$, and rewriting everything in terms of the generating series $B(u,t) := \sum B_k(u)t^k/k!$, $S(u,t) := \sum S_k(u)t^k/k! \in \mathbb{Q}[[u,t]]$, we can write the blowup formulas as

$$\Psi_C^{\widehat{X},M}(xz + te, p^r) = \Psi_C^{X,M}(xz, B(p,t)p^r), \quad \Psi_{C+E}^{\widehat{X},M}(xz + te, p^r) = \Psi_C^{X,M}(xz, S(p,t)p^r). \tag{4.4.2}$$

We now show that these formulas are compatible with our extension of the Donaldson invariants to $\overline{\mathbb{H}}_X$, and give a formula for the power series $B(u,t)$ and $S(u,t)$ in terms of theta functions. [F-S1] also gave explicit formulas for these power series, but in terms of elliptic functions. It is a (not completely trivial) exercise in elliptic functions to show that these formulas are equivalent to ours. However, our formulation, which we derive directly, is more practical for our purposes.

**Proposition 4.5.** *The power series $B(u,t)$, $S(u,t)$ are determined by*

$$B(U(\tau),t) = e^{t^2 G(\tau)/f(\tau)} \theta_{00}(\tau, \tfrac{t}{2\pi i f(\tau)})/\theta(\tau), \tag{4.5.1}$$

$$S(U(\tau),t) = 1^{-1/8} e^{t^2 G(\tau)/f(\tau)} \theta_{11}(\tau, \tfrac{t}{2\pi i f(\tau)})/\theta(\tau). \tag{4.5.2}$$

*Proof.* We determine the formulas by making use of the universality. Let $X$ be a 4-manifold with $b_+ = 1$. We denote by $\widehat{X}$ the connected sum $X \# \bar{\mathbb{P}}_2$ and by $E$ the generator of $H^2(\bar{\mathbb{P}}_2, \mathbb{Z})$, similarly let $\widetilde{X} := \widehat{X} \# \bar{\mathbb{P}}_2$, and let $F$ be the corresponding generator. Let $C \in H^2(X,\mathbb{Z})$, $C \notin 2H^2(X,\mathbb{Z})$. Let $H, M \in \mathbb{H}_X$ with $M\xi \ne 0$ for all $\xi \in H^2(X,\mathbb{Z}) + C/2$. By the proof of Lemma 4.1 of [G], Theorem 4.4 implies that for all $x \in H_2(X,\mathbb{C})$

$$\sum_{H\xi<0<M\xi} (-1)^{\xi C} \delta^X_\xi(e^{xz} p^r B(p,t)) = \sum_{H\xi<0<M\xi} (-1)^{\xi C} \sum_{n \in \mathbb{Z}} \delta^{\widehat{X}}_{\xi+nE}(e^{xz+te} p^r)$$

$$\sum_{H\xi<0<M\xi} (-1)^{\xi C} \delta^X_\xi(e^{xz} p^r S(p,t)) = 1^{-1/8} \sum_{H\xi<0<M\xi} (-1)^{\xi C} \sum_{n \in \mathbb{Z}} (-1)^n \delta^{\widehat{X}}_{\xi+(n+1/2)E}(e^{xz+te} p^r),$$

the sums running as usual over $\xi \in H^2(X,\mathbb{Z}) + C/2$ (note again the different conventions from [G]). Now using Conjecture 2.5 in the same way as in the proof of Lemma 4.5 of [G], we see that we can remove the sums $\sum_{H\xi<0<M\xi}$ on both sides in both equalities, so that

$$\delta^X_\xi(e^{xz} p^r B(p,t)) = \sum_{n \in \mathbb{Z}} \delta^{\widehat{X}}_{\xi+nE}(e^{xz+te} p^r)$$

for all $\xi \in \frac{1}{2} H^2(X,\mathbb{Z})$ with $\xi^2 < 0$, and similarly for $S(p,t)$. Specializing this to $x = 0$ and applying Theorem 2.8, we get

$$\operatorname{Coeff}_{u(\tau)^{r+1}}\left[q^{-\xi^2/2} \theta(\tau)^{\sigma(X)} B(U(\tau),t)/f(\tau)\right]$$

$$= \operatorname{Coeff}_{u(\tau)^{r+1}}\left[q^{-\xi^2/2} \sum_{n \in \mathbb{Z}} q^{n^2/2} e^{nt/f(\tau)} e^{t^2 G(\tau)/f(\tau)^2} \theta(\tau)^{\sigma(X)-1}/f(\tau)\right] \tag{4.5.3}$$

$$= \operatorname{Coeff}_{u(\tau)^{r+1}}\left[q^{-\xi^2/2} \theta_{00}(\tau, \tfrac{t}{2\pi i f(\tau)}) e^{t^2 G(\tau)/f(\tau)^2} \theta(\tau)^{\sigma(X)-1}/f(\tau)\right]$$

If we assume that $\sigma(X) < 0$ (as we may, since the formula is supposed to be universal), then all negative integers appear as $4\xi^2$, so this last formula holds for all $r$ and with $\xi^2$ replaced by $-N/4$ for any integer $N > 0$. As $f(\tau)$ and $u(\tau)$ are power series in $q^{1/2}$ with nonvanishing constant term multiplied with $q^{1/8}$ and $q^{1/4}$, respectively, and $\theta_{00}(\tau,t)$ is even in $t$, this implies the identity (4.5.1). The same argument and the fact that $\theta_{11}(\tau,t)$ is odd in $t$ imply the identity (4.5.2). $\square$

**Proposition 4.6.** *Let $X$ be a simply connected 4-manifold with $b_+ = 1$. Then (4.4.2) holds for all $M \in \overline{\mathbb{H}}_X$ and all $C \in H^2(X,\mathbb{Z})$.*

*Proof.* First we want to remove the assumption that $C$ is not congruent to $0$ modulo $2$. Let $M \in \mathbb{H}_X$ with $M\xi \ne 0$ for all $\xi \in H^2(X,\mathbb{Z})$, fix $d, k \in \mathbb{Z} \ge 0$. By (4.4.1) with $E$ replaced by $F$ together with formula (2.5.2), we have for all $\alpha \in A_d(X)$ and all $\epsilon > 0$ sufficiently small

$$\Phi_0^{X,M}(\alpha B_k(p)) = \Phi_E^{\widehat{X},M+\epsilon E}(e\alpha B_k(p)) = \Phi_E^{\widetilde{X},M+\epsilon E}(e\alpha f^k) = \Phi_0^{\widehat{X},M}(\alpha e^k),$$

$$\Phi_0^{X,M}(\alpha S_k(p)) = \Phi_E^{\widehat{X},M+\epsilon E}(e\alpha S_k(p)) = \Phi_{E+F}^{\widetilde{X},M+\epsilon E}(e\alpha f^k) = \Phi_E^{\widehat{X},M}(\alpha e^k).$$

To show that the blowup formulas hold for all $M \in \overline{\mathbb{H}}_X$ amounts to showing the compatiblity with Definition 4.1. Let $L$ be a unimodular lattice and $h, m \in C_L \cup S_L$. Let $L_1$ be the orthogonal direct sum $L \oplus \langle e \rangle$, where $Q(e) = 1/2$. In view of Proposition 4.5, the compatibility with Definition 4.1 amounts to the identities

$$\Theta^{h,m}_{L_1,c,c}(\tau, xz + te) = \theta_{00}(\tau, t)\, \Theta^{h,m}_{L,c,c}(\tau, xz), \quad \Theta^{h,m}_{L_1,c+e,c+e}(\tau, xz+te) = 1^{1/4}\theta_{11}(\tau,t)\, \Theta^{H,M}_{L,c,c}(\tau, xz)$$

for $c \in L$ and $x \in L$, and these are obvious by the definition of $\Theta^{h,m}_{L,c,c}$ (eq. (3.3.1)). $\square$

**4.3. The structure theorem for the differences.** Let $X$ be a simply connected 4-manifold with $b_+ = 1$ and $F, G \in \mathbb{S}_X$, and use the same letters for their primitive representatives in $H^2(X, \mathbb{Z})$. As we want to describe the differences of the Donaldson invariants at $F$ and at $G$ we write

$$\Phi^{X,F,G}_C := \Phi^{X,F}_C - \Phi^{X,G}_C, \quad \Psi^{X,F,G}_C := \Psi^{X,F}_C - \Psi^{X,G}_C, \quad \overline{\Psi}^{X,F,G}_C := \overline{\Psi}^{X,F}_C - \overline{\Psi}^{X,G}_C.$$

**Definition 4.7.** Let $M = M(X, F, G) := \max\{W^2 \mid W \in H^2(X, \mathbb{Z}) \text{ characteristic}, WF \geq 0 \geq WG\}$. Note that $M \leq 0$, $M \equiv \sigma(X) \pmod{8}$. Put

$$B = B(X, F, G) := \{W \in H^2(X, \mathbb{Z}) \text{ characteristic} \mid W^2 > \sigma(X), (WF)(WG) \leq 0\},$$
$$B_I = B_I(X, F, G) := \{W \in B \mid WF > 0 > WG\},$$
$$B_F = B_F(X, F, G) := \{W \in B \mid 0 = WF < WG \leq 2FG\},$$
$$B_G = B_G(X, F, G) := \{W \in B \mid 0 = WG < WF \leq 2FG\}$$

(cf. Definition 3.12). The elements $W \in B$ are the *basic classes* for $(X, F, G)$. For $(\tau, x) \in \mathfrak{H} \times L_\mathbb{C}$ we put

$$\Omega^{X,F,G}_C(\tau, x) := - \sum_{W \in B_F} \frac{(-1)^{C(W+C)/2} q^{W^2/8} e^{Wx}}{1 - (-1)^{CF} e^{2Fx}} + \sum_{W \in B_G} \frac{(-1)^{C(W+C)/2} q^{W^2/8} e^{Wx}}{1 - (-1)^{CG} e^{2Gx}} \qquad (4.7.1)$$
$$+ \sum_{W \in B_I} q^{W^2/8}\left((-1)^{C(W+C)/2} e^{Wx} - (-1)^{C(-W+C)/2} e^{-Wx}\right).$$

For $x \in L_\mathbb{C}$ we also write

$$O^{X,F,G}_C(x) := - \sum_{\substack{W \in B_F \\ W^2 = M}} \frac{(-1)^{C(W+C)/2} e^{Wx}}{1 - (-1)^{CF} e^{2Fx}} + \sum_{\substack{W \in B_G \\ W^2 = M}} \frac{(-1)^{C(W+C)/2} e^{Wx}}{1 - (-1)^{CG} e^{2Gx}} \qquad (4.7.2)$$
$$+ \sum_{\substack{W \in B_I \\ W^2 = M}} \left((-1)^{C(W+C)/2} e^{Wx} - (-1)^{C(-W+C)/2} e^{Wx}\right).$$

We put

$$A(\tau, x \cdot z) = \theta^0_{10}(\tau)^{\sigma(X)} \frac{4\eta(2\tau)^2}{\eta(\tau)^4} \exp\left(-\mathbf{Q}(x)\big(4G_2(\tau) + 2e_1(\tau)\big)\frac{\eta(2\tau)^4}{\eta(\tau)^8} z^2\right).$$

**Theorem 4.8.** *(Structure theorem)*. *Let $x \in H_2(X, \mathbb{Z})$. If $B_F \neq \emptyset$ and $CF$ is even, assume that $Fx \neq 0$ and if $B_G \neq \emptyset$ and $CG$ is even, assume that $Gx \neq 0$. Let $k = (M - \sigma(X))/8$.*

1. *$\Phi^{X,F,G}_C$ fulfills the $k$-th order simple type condition, i.e. $\Phi^{X,F,G}_C$ vanishes identically if $k \leq 0$ and $\Phi^{X,F,G}_C(\alpha(p^2 - 4)^k) = 0$ for all $\alpha \in A_*(X)$ if $k > 0$.*

2.
$$\overline{\Psi}_C^{X,F,G}(x \cdot z, t) = \sum_{n \geq -1} \left( P_n\left(\frac{t}{1-2t}\right) - 1^{-(C^2+1+n)/4} P_n\left(\frac{-t}{1+2t}\right) \right) z^n$$

for suitable polynomials $P_n(y)$ $(n \geq -1)$ of degree $\leq k$ with no constant term.

3. The polynomials $P_n(y)$ are determined by the fact that $\overline{\Psi}_C^{X,F,G}(x \cdot z, \widetilde{u}(\tau))$ is the principal part in the development of $A(\tau, x \cdot z) \Omega_C^{X,F,G}\left(\tau, xz \frac{\eta(2\tau)^2}{\eta(\tau)^4}\right)$ as a Laurent series in $\widetilde{U}(\tau) - 2$.

4. $\Psi_C^{X,F,G}(x \cdot z, (1+p/2)(p^2-4)^{k-1}) = 2^{1+M} O_C^{X,F,G}(x \cdot z) e^{\mathbf{Q}(x)z^2/2}$.

*Proof.* We apply Theorem 3.13 to $\varphi_{L,c}^{f,g}(\tau, x \cdot z)$ for $L = H_2(X, \mathbb{Z})$ with the negative of the intersection form and $f$, $g$, $c$ the Poincaré duals of $F$, $G$ and $C$, respectively. We obtain

$$\overline{\Psi}_C^{X,F,G}(x \cdot z, u(\tau)) = \varphi_{L,c}^{f,g}(\tau, x \cdot z) - \sum_{n \geq -1} R_n(U(\tau)) z^n,$$

where $R_n(t)$ is the polynomial from Theorem 3.13. Therefore 2. and 3. follow directly from parts 2. and 3. of Theorem 3.13. To show 1. note that by 2.

$$\Phi_C^{X,F,G}\left(p^r(p^2-4)^k x^n/n!\right) = \text{Coeff}_{z^n t^{r+1}}\left[(t^{-2}-4)^k \overline{\Psi}_C^{X,F,G}(x \cdot z, t)\right]$$
$$= \text{Coeff}_{t^{r+1}}\left[(t^{-2}-4)^k \left(P_n\left(\frac{1}{t^{-1}-2}\right) - 1^{-(C^2+1+n)/4} P_n\left(\frac{-1}{t^{-1}+2}\right)\right)\right],$$

which is 0 because $P_n(y)$ has degree $\leq k$ and has no constant term.

To prove 4. let $a_n$ be the coefficient of the leading term in $P_n(y)$ as in Theorem 3.13. Then we get

$$\Psi_C^{X,F,G}(x \cdot z, (1+p/2)(p^2-4)^{k-1}) = \frac{1}{2}\text{Coeff}_t\Bigg[(t^{-1}+2)^k (t^{-1}-2)^{k-1}$$
$$\cdot \sum_{n \geq -1} \left(P_n\left(\frac{1}{t^{-1}-2}\right) - 1^{-(C^2+1+n)/4} P_n\left(\frac{-1}{t^{-1}+2}\right)\right)\Bigg]$$
$$= \frac{1}{2}\text{Coeff}_t\left[(t^{-1}+2)^k \sum_{n \geq -1} a_n (t^{-1}-2)^{-1} z^n\right] = 2^{2k-1} \sum_{n \geq -1} a_n z^n.$$

So the result follows by part 4. of Theorem 3.13. $\square$

**Remark 4.9.** 1. Now we view $z$ as a complex variable. By the obvious identity

$$\frac{(-1)^{C(C+W)/2} e^{Wxz}}{1 - (-1)^{CF} e^{2Fxz}} = \sum_{n \geq 0} (-1)^{C(C+W+2nF)/2} e^{(W+2nF)xz} \qquad (\Re(Fxz) < 0)$$

we can (again for $\Re(Fxz) < 0$) write

$$\Omega_C^{X,F,G}(\tau, x \cdot z) = \sum_{W \in B} (-1)^{C(C+W)/2} \epsilon_W q^{W^2/8} e^{Wxz}, \qquad \epsilon_W = \begin{cases} 1 & WF > 0 \geq WG \\ -1 & WF \leq 0 < WG \\ 0 & \text{otherwise.} \end{cases}$$

Therefore we get that $\overline{\Psi}_C^{X,F,G}(x \cdot z, \widetilde{u}(\tau))$ is the principal part of the development of

$$\sum_{W \in B} (-1)^{C(C+W)/2} \epsilon_W \, q^{W^2/8} \exp\left(Wxz\frac{\eta(2\tau)^2}{\eta(\tau)^4}\right) A(\tau, x \cdot z)$$

as a Laurent series in $\widetilde{U}(\tau) - 2$, when the latter converges. So we can view the elements of $B$ as basic classes in the sense of the Introduction. Note that in general the signs $\epsilon_W$ will depend on the signs of $\Re(Fxz)$, $\Re(Gxz)$

2. The principal part of the Laurent development of $\exp\left(Wxz\frac{\eta(2\tau)^2}{\eta(\tau)^4}\right)A(\tau, x \cdot z)$ (i.e., up to a factor $(-1)^{C(C+W)/2}\epsilon_W$, the contribution of the basic class $W$ to $\overline{\Psi}_C^{X,F,G}(x \cdot z, \widetilde{u}(\tau))$) is

$$e^{Wxz} e^{\mathbf{Q}(x)z^2/2} Q_{k_W}((\widetilde{U}(\tau) - 2)^{-1}, \sigma(X), z^2\mathbf{Q}(x), Wxz)$$

where $k_W = (W^2 - \sigma(X))/8$ and $Q_k(z_1, z_2, z_3, z_4)$ is a universal polynomial of degree $k$. This can be seen by writing

$$A(\tau, x \cdot z) = q^{(\sigma(X)-W^2)/8}(q^{-1/8}\theta_{10}^0(\tau))^{\sigma(X)} \frac{4\eta(2\tau)^2}{\eta(\tau)^4} \exp\left(-\mathbf{Q}(x)(4G_2(\tau) + 2e_1(\tau))\frac{\eta(2\tau)^4}{\eta(\tau)^8}z^2\right),$$

and using the developments

$$q = \frac{y}{64} - \frac{y^2}{512} + \ldots, \qquad q^{-1/8}\theta_{10}^0(\tau) = 2 + \frac{y}{32} - \frac{y^2}{256} + \ldots,$$

$$\frac{\eta(2\tau)^2}{\eta(\tau)^4} = 1 + \frac{y}{16} - \frac{5y^2}{1024} + \ldots, \quad (4G_2(\tau) + 2e_1(\tau))\frac{\eta(2\tau)^4}{\eta(\tau)^8} = \frac{1}{2} + \frac{y}{8} - \frac{y^2}{256} + \ldots.$$

in $y := \widetilde{U}(\tau) - 2$. The only term in the expression for $A(\tau, x \cdot z)$ containing negative powers of $y$ will then be $q^{(\sigma(X)-W^2)/8}$. Writing

$$\exp\left(bz\left(1 + \sum a_i y^i\right)\right) = e^{bz} \cdot \exp\left(bz\left(\sum a_i y^i\right)\right)$$

and expanding the second factor, the result follows. Using part 2. of Theorem 4.8 we can also rewrite this result as follows: Up to a factor $(-1)^{C(C+W)/2}\epsilon_W$ the contribution of the basic class $W$ to $\overline{\Psi}_C^{X,F,G}(x \cdot z, t)$ is

$$e^{Wxz} e^{\mathbf{Q}(x)/2} Q_{k_W}((t^{-1} - 2)^{-1}, \sigma(X), z^2\mathbf{Q}(x), Wxz)$$
$$- 1^{-C^2-1} e^{-iWxz} e^{-\mathbf{Q}(x)/2} Q_{k_W}((-t^{-1} - 2)^{-1}, \sigma(X), -z^2\mathbf{Q}(x), -iWxz).$$

3. As we noted in Remark 2.7, $\Gamma_u$ is conjugate to $\Gamma(2)$ via a matrix $M \in GL(2, \mathbb{Z})$ sending $\infty, 1, -1$ to $1, 0, \infty$ respectively. Therefore we could reexpress our results in terms of developments of modular functions for $\Gamma(2)$ in powers of the modular function $\bar{u}(\tau) = u(M\tau)$, thus getting a connection to the computations [W1], [W2] in theoretical physics.

**Corollary 4.10.** *If $\sigma(X) > -8$, then $\overline{\Psi}_C^{X,F}$ is independent of $F \in \mathbb{S}_X$.*

*Proof.* This is immediate from part 2. of Theorem 4.8, since $k < 1$. $\square$

We also get another version of the blowup formulas.

**Corollary 4.11.** *Let $F, G \in \mathbb{S}_X$ and $C \in H^2(X, \mathbb{Z})$ and $x \in H_2(X, \mathbb{Z})$; we denote by the same letters the pullbacks to $\widehat{X} = X \# \overline{\mathbb{P}}_2$, and by $E$ the class of the exceptional divisor. Let $t$ be an indeterminate and set $s = t\, \eta(2\tau)^2/\eta(\tau)^4$. Then we have*

1. *With respect to the developments in powers of $\widetilde{U}(\tau) - 2$, the function $\overline{\Psi}_C^{\widehat{X}, F, G}(x \cdot z + te, \widetilde{u}(\tau))$ is the principal part of*

$$\overline{\Psi}_C^{X, F, G}\bigl(x \cdot z, \widetilde{u}(\tau)\bigr) \, \exp\bigl(s^2(4G_2(\tau) + 2e_1(\tau))\bigr) \cosh(s) \prod_{n > 0} \frac{(1 + q^n e^{2s})(1 + q^n e^{-2s})}{(1 + q^n)^2}$$

*and $\overline{\Psi}_{C+E}^{\widehat{X}, F, G}(x \cdot z + te, \widetilde{u}(\tau))$ is the principal part of*

$$\overline{\Psi}_C^{X, F, G}\bigl(x \cdot z, \widetilde{u}(\tau)\bigr) \, \exp\bigl(s^2(4G_2(\tau) + 2e_1(\tau))\bigr) \sinh(s) \prod_{n > 0} \frac{(1 - q^n e^{2s})(1 - q^n e^{-2s})}{(1 + q^n)^2}.$$

2. *In particular, putting $M := M(X, F, G)$ and $k := (M - \sigma(X)/8)$, we get $\Phi_D^{\widehat{X}, F, G}(\alpha(p^2 - 4)^k) = 0$ for all $\alpha \in A_*(\widehat{X})$ and all $D \in H^2(\widehat{X}, \mathbb{Z})$ and*

$$\Phi_C^{\widehat{X}, F, G}\bigl(e^{xz + te}(1 + p/2)(1 - p^2/4)^{k-1}\bigr) = \Phi_C^{X, F, G}\bigl(e^{xz}(1 + p/2)(1 - p^2/4)^{k-1}\bigr) \cosh(t)\, e^{-t^2/2},$$

$$\Phi_{C+E}^{\widehat{X}, F, G}\bigl(e^{xz + te}(1 + p/2)(1 - p^2/4)^{k-1}\bigr) = \Phi_C^{X, F, G}\bigl(e^{xz}(1 + p/2)(1 - p^2/4)^{k-1}\bigr) \sinh(t)\, e^{-t^2/2}.$$

*Proof.* This follows from Propositions 4.5 and 4.6 by applying $W$ and using the standard identities

$$\bigl(\theta^{-1}\theta_{00} | W\bigr)(\tau, t) = \cosh(\pi i t) \prod_{n > 0} \frac{(1 + q^n e^{2\pi i t})(1 + q^n e^{-2\pi i t})}{(1 + q^n)^2},$$

$$\bigl(\theta^{-1}\theta_{11} | W\bigr)(\tau, t) = 1^{1/8} \sinh(\pi i t) \prod_{n > 0} \frac{(1 - q^n e^{2\pi i t})(1 - q^n e^{-2\pi i t})}{(1 + q^n)^2}.$$

We omit the details, which are not difficult. Notice that the principle used here (as already in the proof of Theorem 3.13) is that if the expansions of two functions at $\infty$ coincide and if both functions are known to be modular, then their expansions at any other cusp also coincide. □

**4.4. Speculations about the relation to Seiberg-Witten theory.** Theorem 4.8 is closely related to the expectations from Seiberg-Witten theory. The elements $W \in B$ are characteristic elements of $H^2(X, \mathbb{Z})$ and thus correspond to $Spin^c$-structures, and $(W^2 - \sigma(X))/8 - 1$ is the expected dimension of the corresponding Seiberg-Witten moduli space. The $SW$-basic classes are those classes $W$, for which the Seiberg-Witten invariant is not zero. $X$ is of $SW$-simple type if only classes for which the expected dimension is 0 give rise to nonzero invariants. At least in the case $b_+ > 1$ the set of $SW$-basic classes together with the corresponding invariants and expected dimensions are conjectured [W1] to determine the Donaldson invariants of $X$ by a universal formula, similar to Theorem 4.8. In case $b_+ > 1$ there is only a finite number of $SW$-basic classes.

In the case of $b_+ = 1$ this relationship should still be there, but obscured by the chamber structure both in Seiberg-Witten and Donaldson theory. According to the program of [P-T1], [P-T2] the Donaldson invariants should be determined from the Seiberg-Witten invariants after a large perturbation of the equation; the Donaldson invariants $\Phi_{C, d}^{X, g}$ for a metric should be more or less given in terms of Seiberg-Witten invariants for perturbed equations, with the perturbation depending on $d$. This also makes it possible that an infinite number of basic classes contribute to $\Phi_C^{X,g}$. Our formulas suggest however that for period points $F, G \in \mathbb{S}_X$, the situation again becomes easier. Denote by $SW_L(W)$ the unperturbed Seiberg-Witten invariant for the $Spin^c$-structure with first Chern-class $W$, for the metric with period point $L \in \mathbb{H}_X$. It is well known that $SW_L(W)$ is constant on both connected components of $\mathbb{H}_X \setminus W^\perp$, and that for $L_+ W > 0 > L_- W$ we have $SW_{L_+}(W) = SW_{L_-}(W) \pm 1$. This also shows that 4-manifolds with $b_+ = 1$ are, for essentially trivial reasons, usually not of $SW$-simple type. Theorem 4.8 says that for each class $W$ for which $SW_L(W)$ changes when going from $F$ to $G$ we get a contribution to the difference $\overline{\Psi}_C^{X,F,G}(x \cdot z, t)$ given by a universal formula. Classes orthogonal to $F$ and $G$ are treated in a special way. Therefore we conjecture that our formula indeed gives the formula for the conjectured relation between Donaldson and Seiberg-Witten invariants, viz.:

**Conjecture 4.12.** *Let $X$ a simply-connected 4-manifold with $b_+ = 1$, $F$ a primitive representative of a class in $\mathbb{S}_X$, and $B_1 = \{W \in H^2(X, \mathbb{Z}) | \ W \ characteristic, \ W^2 > \sigma(X), \ WF = 0\}$. Then there is a system $R_F$ of representatives of $B_1$ modulo $\langle 2F \rangle$ such that $\Psi_C^{X,F}(\widetilde{u}(\tau), x \cdot z)$ is the principal part in the development of $A(\tau, xz) \Omega_C^{X,F}\left(\tau, xz \frac{\eta(2\tau)^2}{\eta(\tau)^4}\right)$ as a Laurent series in $\widetilde{U}(\tau) - 2$, where*

$$\Omega_C^{X,F}(\tau, x) := - \sum_{W \in R_F} \frac{(-1)^{C(W+C)/2} q^{W^2/8} e^{Wx}}{1 - (-1)^{CF} e^{2Fx}} + \sum_W q^{W^2/8} (-1)^{C(W+C)/2} SW_{N_W}(W) e^{Wx}.$$

*In the second sum $W$ runs through the characteristic elements of $H^2(X, \mathbb{Z})$ and for each $W$ the period point $N_W \in \mathbb{H}_X$ is chosen in such a way that $(WN_W)(WF) > 0$ (for a representative of $N_W$ with $N_W F > 0$). If $CF$ is odd and $Fx = 0$, we can replace the first sum by $-\sum_W (-1)^{C(W+C)/2} q^{W^2/8} e^{Wx}/2$, with $W$ now running through any system of representatives of $B_1$ modulo $\langle 2F \rangle$.*

Note that this is compatible with the predictions of Witten in the simple type case. One can suspect that a modification of this formula should work in the case $b_+ > 1$ (if 4-manifolds with $b_+ > 1$ not of simple type do indeed exist).

## 5. The case of rational algebraic surfaces

For rational algebraic surfaces $X$ we shall see that there are always some $G \in \mathbb{S}(X)$ such that $\Phi_C^{X,G} = 0$ for all $C \in H^2(X, \mathbb{Z})$. Therefore Theorem 4.8 will give us the structure of the Donaldson invariants $\Phi_C^{X,F}$ for $F \in \mathbb{S}_X$ instead of only the differences.

**5.1. Donaldson invariants of $\mathbb{P}_1 \times \mathbb{P}_1$ and $\mathbb{P}_2 \# \overline{\mathbb{P}}_2$.** As a first application of Theorem 4.8 we want to study the Donaldson invariants of $\mathbb{P}_1 \times \mathbb{P}_1$ and $\mathbb{P}_2 \# \overline{\mathbb{P}}_2$, and thus of all rational ruled surfaces. We compute the limit of these Donaldson invariants for the period point going to the boundary of the positive cone. We also show that they satisfy certain relations, also for period points in the inside of the positive cone. We use the following elementary result from algebraic geometry:

**Lemma 5.1.** *[Q2] Let $X \longrightarrow \mathbb{P}_1$ be a rational ruled surface, $F \in H^2(X, Z)$ the class of a fibre and $H$ the class of a section of the ruling. If $C \in H^2(X, \mathbb{Z})$ fulfills $CF$ odd, then the moduli space $M_{NF+H}(C, c_2)$ of $NF + H$ stable rank 2 sheaves on $X$ is empty for all $N \in \mathbb{Z}_{>0}$ sufficiently large with respect to the second Chern class $c_2$. In particular, for any given $d$ the invariant $\Phi_{C,d}^{X,F+\epsilon H}$ vanishes for all sufficiently small $\epsilon > 0$.*

**Notation 5.2.** 1. Let $F$ and $G$ be the Poincaré duals of the classes of the fibres of the projections of $\mathbb{P}_1 \times \mathbb{P}_1$ to its factors. We denote by $\widehat{\mathbb{P}}_2$ the blowup of $\mathbb{P}_2$ in a point and by $E_1$ the class of the exceptional divisor. We denote by $\bar{F} := H - E_1$ the Poincaré dual of the fibre of the ruling of $\widehat{\mathbb{P}}_2$ and $\bar{G} := H + E_1$, for $H$ the pullback of the hyperplane class. Let $\sigma : \widetilde{\mathbb{P}}_2 \to \widehat{\mathbb{P}}_2$ be the blowup in a general point, with exceptional divisor $E_2$. Then there is a blowup $\epsilon : \widetilde{\mathbb{P}}_2 \to \mathbb{P}_1 \times \mathbb{P}_1$ in a point with exceptional divisor $E$ such that $\epsilon^*(F) = \sigma^*(\bar{F})$, $\epsilon^*(G) = \sigma^*(H) - E_2$, $E = \epsilon^*(F) - E_2$.

2. We denote by $f, g, \bar{f}, \bar{g}$ the Poincaré duals of $F, G, \bar{F}, \bar{G}$ respectively.

3. For $X = \mathbb{P}_1 \times \mathbb{P}_1$ or $X = \widehat{\mathbb{P}}_2$ and $L \in \{F, G, \bar{F}, \bar{G}\}$, $C \in H^2(X, \mathbb{Z})$ and $\alpha \in A_d(X)$ we denote by $\Phi^{X, L+}_{C, d}(\alpha) := \Phi^{X, L+\epsilon M}_{C, d}(\alpha)$ for $M$ an ample divisor and $\epsilon > 0$ sufficiently small. As before we put $\Phi^{X, L+}_C := \sum_d \Phi^{X, L+}_{C, d}$. Note that if $LC$ is odd then $\Phi^{X, L+}_C = \Phi^{X, L}_C$.

**Theorem 5.3.** 1. *For $X = \mathbb{P}_1 \times \mathbb{P}_1$ and $L \in \{F, G\}$ or $X = \widehat{\mathbb{P}}_2$, and $L \in \{\bar{F}, \bar{G}\}$ and for all $C \in H^2(X, \mathbb{Z})$, the Donaldson invariants $\Phi^{X, L}_C$ vanish.*

2. *For all $r \geq 0$ and indeterminates $s, t$ we have*

$$\Phi^{\mathbb{P}_1 \times \mathbb{P}_1, F+}_0(e^{sf+tg}p^r) = -\operatorname{Coeff}_{u(\tau)^{r+1}}\left[\coth\left(\tfrac{t}{2f(\tau)}\right)e^{-2stG(\tau)/f(\tau)^2}/f(\tau)\right],$$

$$\Phi^{\widehat{\mathbb{P}}_2, \bar{F}+}_0(e^{s\bar{f}+t\bar{g}}p^r) = \Phi^{\mathbb{P}_1 \times \mathbb{P}_1, F+}_0(e^{sf+2tg}p^r),$$

$$\Phi^{\mathbb{P}_1 \times \mathbb{P}_1, F+}_F(e^{sf+tg}p^r) = -\operatorname{Coeff}_{u(\tau)^{r+1}}\left[\frac{e^{-2stG(\tau)/f(\tau)^2}}{\sinh\left(\tfrac{t}{2f(\tau)}\right)f(\tau)}\right],$$

$$\Phi^{\widehat{\mathbb{P}}_2, \bar{F}+}_{\bar{F}}(e^{sf+tg}p^r) = \Phi^{\mathbb{P}_1 \times \mathbb{P}_1, F+}_0(e^{sf+2tg}p^r).$$

3. *More generally we have for all $a, b \in \mathbb{R}_{>0}$*

$$\Phi^{\mathbb{P}_1 \times \mathbb{P}_1, aF+bG}_F(e^{sf+tg}p^r) = \Phi^{\mathbb{P}_1 \times \mathbb{P}_1, aG+bF}_G(e^{sg+tf}p^r),$$

$$\Phi^{\mathbb{P}_1 \times \mathbb{P}_1, aF+bG}_0(e^{sf+tg}p^r) = \Phi^{\mathbb{P}_1 \times \mathbb{P}_1, aG+bF}_0(e^{sg+tf}p^r),$$

$$\Phi^{\widehat{\mathbb{P}}_2, a\bar{F}+bG}_{\bar{F}}(e^{s\bar{f}+t\bar{g}}p^r) = \Phi^{\mathbb{P}_1 \times \mathbb{P}_1, aF+2bG}_F(e^{sf+2tg}p^r) - \Phi^{\mathbb{P}_1 \times \mathbb{P}_1, 2aF+bG}_G(e^{2sf+tg}p^r)$$

$$= \Phi^{\mathbb{P}_1 \times \mathbb{P}_1, 2aF+bG}_0(e^{2sf+tg}p^r) - \Phi^{\mathbb{P}_1 \times \mathbb{P}_1, aF+2bG}_0(e^{sf+2tg}p^r),$$

$$\Phi^{\widehat{\mathbb{P}}_2, a\bar{F}+b\bar{G}}_0(e^{s\bar{f}+t\bar{g}}p^r) = \Phi^{\mathbb{P}_1 \times \mathbb{P}_1, 2aF+bG}_0(e^{2sf+tg}p^r) - \Phi^{\mathbb{P}_1 \times \mathbb{P}_1, aF+2bG}_F(e^{sf+2tg}p^r).$$

One can check that these formulas are compatible up to conventions with the known results (e.g. [L-Q],[K-L],[E-G2]). Note that part 2 of the theorem implies in particular Conjecture 6.2 from [E-G2].

*Proof.* 1. By Lemma 5.1 we get that

$$\Phi^{\mathbb{P}_1 \times \mathbb{P}_1, F}_{F+G} = \Phi^{\mathbb{P}_1 \times \mathbb{P}_1, F}_G = \Phi^{\mathbb{P}_1 \times \mathbb{P}_1, G}_F = \Phi^{\widehat{\mathbb{P}}_2, \bar{F}}_H = \Phi^{\widehat{\mathbb{P}}_2, \bar{G}}_H = \Phi^{\widehat{\mathbb{P}}_2, \bar{F}}_{E_1} = \Phi^{\widehat{\mathbb{P}}_2, \bar{G}}_{E_1} = 0.$$

Applying Corollary 4.10 it follows that $\Phi^{\mathbb{P}_1 \times \mathbb{P}_1, F}_F = \Phi^{\mathbb{P}_1 \times \mathbb{P}_1, G}_G = 0$. In the following let $x \in H_2(\mathbb{P}_1 \times \mathbb{P}_1, \mathbb{C})$ and $y \in H_2(\widehat{\mathbb{P}}_2, \mathbb{C})$. We apply Proposition 4.6 to obtain

$$\Phi^{\widetilde{\mathbb{P}}_2, \epsilon^*(F)}_{\epsilon^*(F)}(e^{xz+te}p^r) = \Phi^{\mathbb{P}_1 \times \mathbb{P}_1, F}_F(e^{xz}B(p, t)p^r) = 0,$$

i.e. $\Phi^{\widetilde{\mathbb{P}}_2, \epsilon^*(F)}_{\epsilon^*(F)} = 0$. Therefore we have by again using 4.6

$$\Phi^{\widetilde{\mathbb{P}}_2, \bar{F}}_{\bar{F}}(e^{yz}p^r) = \Phi^{\widetilde{\mathbb{P}}_2, \epsilon^*(F)}_{\epsilon^*(F)}(e^{yz}p^r) = 0.$$

Furthermore we have by repeating essentially the same argument

$$\Phi^{\widetilde{\mathbb{P}}_2, \epsilon^*(F)}_{E_2}(e^{xz+te}p^r) = -\Phi^{\mathbb{P}_1 \times \mathbb{P}_1, F}_{F}(e^{xz}S(p,t)p^r) = 0,$$

and thus $\Phi^{\widehat{\mathbb{P}}_2, \bar{F}}_0(e^{yz}p^r) = \Phi^{\widetilde{\mathbb{P}}_2, \epsilon^*(F)}_{E_2}(e_2 e^{yz}p^r) = 0$. Finally we get

$$\Phi^{\widetilde{\mathbb{P}}_2, \epsilon^*(F)}_E(e^{yz+te}p^r) = -\Phi^{\widehat{\mathbb{P}}_2, F}_F(e^{yz}S(p,t)p^r) = 0,$$

and thus $\Phi^{\mathbb{P}_1 \times \mathbb{P}_1, F}_0(e^{xz}p^r) = \Phi^{\widetilde{\mathbb{P}}_2, \epsilon^*(F)}_E(e^{xz}ep^r) = 0$. The other cases follow by symmetry.

2. Putting $X = \mathbb{P}_1 \times \mathbb{P}_1$ and $L := aF + bG$ for $a, b \in \mathbb{R}_{>0}$ and writing $\nu(t) := \begin{cases} 1 & t > 0, \\ 1/2 & t = 0 \\ 0 & t < 0, \end{cases}$ we get by part 1., Definition 4.1 and formula (3.3.2) the formulas

$$\Phi^{\mathbb{P}_1 \times \mathbb{P}_1, L}_0(e^{xz}p^r) = \sum_{n,m>0} \nu(bn - am)\, \delta^{\mathbb{P}_1 \times \mathbb{P}_1}_{nF-mG,d}(e^{xz}p^r) - \mathrm{Coeff}_{u(\tau)^{r+1}}[L_0(\tau, xz)],$$

$$\Phi^{\mathbb{P}_1 \times \mathbb{P}_1, L}_F(e^{xz}p^r) = \sum_{n,m>0} (-1)^m \nu\!\left(b(n - \tfrac{1}{2}) - am\right) \delta^{\mathbb{P}_1 \times \mathbb{P}_1}_{(n-1/2)F-mG,d}(e^{xz}p^r) - \mathrm{Coeff}_{u(\tau)^{r+1}}[L_F(\tau, xz)],$$

where we have put

$$L_0(\tau, xz) = e^{-\mathbf{Q}(x)z^2 G(\tau)/f(\tau)^2} \frac{1 + e^{-Fxz/f(\tau)}}{f(\tau)\left(1 - e^{-Fxz/f(\tau)}\right)},$$

$$L_F(\tau, xz) = 2 \frac{e^{-\mathbf{Q}(x)z^2 G(\tau)/f(\tau)^2}}{f(\tau)\left(e^{Fxz/(2f(\tau))} - e^{-Fxz/(2f(\tau))}\right)}.$$

If $a/b$ is sufficiently large, then for all $n, m$ ocurring in the sums above the number $nm$ is larger then $(d+3)/4$ and thus $\delta^{\mathbb{P}_1 \times \mathbb{P}_1}_{nF-mG,d} = 0$ and $\delta^{\mathbb{P}_1 \times \mathbb{P}_1}_{(n-1/2)F-mG,d} = 0$. Therefore

$$\Phi^{\mathbb{P}_1 \times \mathbb{P}_1, F+}_0(e^{xz}p^r) = -\mathrm{Coeff}_{u(\tau)^{r+1}}[L_0(\tau, xz)], \quad \Phi^{\mathbb{P}_1 \times \mathbb{P}_1, F+}_F(e^{xz}p^r) = -\mathrm{Coeff}_{u(\tau)^{r+1}}[L_F(\tau, xz)].$$

The argument for $\widehat{\mathbb{P}}_2$ is analogous.

3. The first two identities are obvious by symmetry. Part 2. implies that, given $d$, the other formulas hold until degree $d$ in $s,t$ for $a/b$ sufficiently large. Note that by Theorem 2.4

$$\Phi_0^{\mathbb{P}_1\times\mathbb{P}_1,aF+bG} - \Phi_0^{\mathbb{P}_1\times\mathbb{P}_1,F+} = \sum_{n,m>0} \nu(bn-am)\delta_{nF-mG}^{\mathbb{P}_1\times\mathbb{P}_1},$$

$$\Phi_F^{\mathbb{P}_1\times\mathbb{P}_1,aF+bG} - \Phi_F^{\mathbb{P}_1\times\mathbb{P}_1,F+} = \sum_{n,m>0} (-1)^m \nu\left(b\left(n-\tfrac{1}{2}\right)-am\right)\delta_{(n-1/2)F-mG}^{\mathbb{P}_1\times\mathbb{P}_1},$$

$$\Phi_0^{\widehat{\mathbb{P}}_2,a\bar{F}+b\bar{G}} - \Phi_0^{\widehat{\mathbb{P}}_2,\bar{F}+} = \sum_{n,m>0} \Big(\nu(bn-am)\delta_{n\bar{F}-m\bar{G}}^{\widehat{\mathbb{P}}_2}$$
$$+ \nu\left(b\left(n-\tfrac{1}{2}\right)-a\left(m-\tfrac{1}{2}\right)\right)\delta_{(n-1/2)\bar{F}-(m-1/2)\bar{G}}^{\widehat{\mathbb{P}}_2}\Big),$$

$$\Phi_{\bar{F}}^{\widehat{\mathbb{P}}_2,a\bar{F}+b\bar{G}} - \Phi_{\bar{F}}^{\widehat{\mathbb{P}}_2,\bar{F}+} = \sum_{n,m>0} \Big(\nu\left(b\left(n-\tfrac{1}{2}\right)-am\right)\delta_{(n-1/2)\bar{F}-m\bar{G}}^{\widehat{\mathbb{P}}_2}$$
$$- \nu\left(bn-a\left(m-\tfrac{1}{2}\right)\right)\delta_{n\bar{F}-(m-1/2)\bar{G}}^{\widehat{\mathbb{P}}_2}\Big).$$

Using Theorem 2.8, this reduces the proof to the easy identities

$$\sum_{n,m>0} \left(\nu\left(b\left(n-\tfrac{1}{2}\right)-am\right)q^{(2n-1)m}e^{(2n-1)t-2ms} - \nu\left(bn-a\left(m-\tfrac{1}{2}\right)\right)q^{n(2m-1)}e^{2nt-(2m-1)s}\right)$$

$$= \sum_{n,m>0} \Big((-1)^m \nu(b(2n-1)-am)q^{(n-1/2)m}e^{(2n-1)t-ms} -$$
$$(-1)^m \nu(bn-a(2m-1))q^{(n-1/2)m}e^{nt-(2m-1)s}\Big)$$

$$= \sum_{n,m>0} \left(\nu(bn-2am)q^{nm}e^{nt-2ms} - \nu(2bn-am)q^{nm}e^{2nt-ms}\right),$$

$$\sum_{n,m>0} \left(\nu(bn-am)q^{2nm}e^{2nt-2ms} + \nu\left(b\left(n-\tfrac{1}{2}\right)-a\left(m-\tfrac{1}{2}\right)\right)q^{(n-1/2)(2m-1)}e^{(2n-1)t-(2m-1)s}\right)$$

$$= \sum_{n,m>0} \left(\nu(bn-2am)q^{nm}e^{nt-2ms} - (-1)^m \nu(b(2n-1)-am)q^{(n-1/2)m}e^{(2n-1)t-ms}\right).$$

$\square$

5.2. **The structure theorem for rational surfaces.** We want to determine the structure of the Donaldson invariants of rational algebraic surfaces $X$ at period points $F$ in the boundary of the positive cone. We already know that all the Donaldson invariants $\Phi_C^{\mathbb{P}_1\times\mathbb{P}_1,F}$ for $F \in \mathbb{S}_{\mathbb{P}_1\times\mathbb{P}_1}$ vanish. As the Donaldson invariants depend only on the diffeomorphism type of the pair $(X,F)$, we can assume that $X$ is $\mathbb{P}_2$ blown up in $N$ points.

**Notation 5.4.** Let $X$ be the blowup of $\mathbb{P}_2$ in $N$ points. Let $H \in H^2(X,\mathbb{Z})$ be the pullback of the hyperplane class from $\mathbb{P}_2$, and denote by $E_1,\ldots,E_N$ the classes of the exceptional divisors. A class $D \in H^2(X,\mathbb{Z})$ is written as

$$D = (d_0,d_1,\ldots,d_N) := d_0 H + \sum_{i=1}^N d_i E_i$$

(uppercase and lowercase letters correspond). If $d_0 \geq d_1 \geq \ldots \geq d_N$, we write also $D = (\alpha_0^{n_0},\ldots,\alpha_k^{n_k})$, where the $\alpha_j$ are the different integers occuring in $d_0,d_1,\ldots,d_N$ with multiplicity $n_j > 0$.

**Corollary 5.5.** 1. $\Phi_C^{X,H+E_j} = 0$ for all $j \in \{1 \ldots N\}$ and all $C \in H^2(X,\mathbb{Z})$.
  2. If $N < 9$, then $\Phi_C^{X,F} = 0$ for all $F \in \mathbb{S}_X$.

*Proof.* 1. holds by Theorem 5.3 in case $N = 1$ and follows by Proposition 4.6 in the general case.
2. follows from 1. by Corollary 4.10. □

**Definition 5.6.** Let $F \in H^2(X,\mathbb{Z})$ be a representative of an element in $\mathbb{S}_X$, with $f_0 \geq f_1$. Let

$$B(X,F) := \{W \in H^2(X,\mathbb{Z}) \mid \text{all } w_i \text{ odd}, W^2 > -N+1, (w_0 - w_1)WF \leq 0\},$$
$$B_I(X,F) := \{W \in B(X,F) \mid w_0 - w_1 < 0 < WF\},$$
$$B_F(X,F) := \{W \in B(X,F) \mid 0 = WF < w_0 - w_1 < 2(f_0 - f_1)\}, B_G(X,F) = \emptyset.$$

The nonnegative integer $(W^2 + N - 1)/8$ is the order of $W \in B(X,F)$. Usually we will drop $(X,F)$ in the notation. We define $\Omega_C^{X,F}(\tau,x)$ and $O_C^{X,F}(x)$ by the formulas (4.7.1) resp. (4.7.2) replacing $B_I(X,F,G)$, $B_F(X,F,G)$, $B_G(X,F,G)$, by $B_I(X,F)$, $B_F(X,F)$, $B_G(X,F)$ respectively.

**Theorem 5.7.** *For a rational surface $X$ and $F \in \mathbb{S}_X$, we can replace $(X,F,G)$ by $(X,F)$ everywhere in Theorem 4.8, Remark 4.9 and Corollary 4.11.*

*Proof.* We put $G := H + E_1$ in Theorem 4.8 and apply Corollary 5.5. □

**Corollary 5.8.** *Let $X$ be a rational surface. Then Conjecture 4.12 holds with $R_F = B_F$. The set of characteristic elements $W \in H^2(X,\mathbb{Z})$ with $SW_{N_W}(W) \neq 0$ is $B_I \cup -B_I$.*

**5.3. Examples.** We want to finish the paper by giving a number of examples which illustrate our main results. Also they should make it clear that for any given pair $(X,F)$ of a rational algebraic surface $X$ and $F \in \mathbb{S}_X$ it is an elementary task to determine all the basic classes for $(X,F)$ and their orders, and thus to compute the corresponding Donaldson invariants completely. We will always give the elements in $B_F$ and $-B_I$.

**Remark 5.9.** The following observations simplify the task of determining the basic classes $W$ for a given $(X,F)$ with $F \in \mathbb{S}_X$.

1. We can always assume that $f_0 \geq f_1 \geq \ldots \geq f_N \geq 0$ (if necessary change the numbering of the $E_i$), and in this case an element $W$ in $B_F \cup -B_I$ must satisfy $w_0 > 0$ and $w_0 > w_i$ for all $i > 0$.
2. If $N - 1 \geq 8$ is divisible by 8 and $(X,F)$ is of strictly $(N-1)/8$-th order simple type, then all the $f_i$ are odd, and $F$ is the only class in $B_F \cup -B_I$ of maximal order: any other such basic class $V$ would have to satisfy $V^2 = 0$ and this already determines $F$.
3. For $W \in B_F$ also $2F - W \in B_F$, and thus we can get all elements $W \in B_F$ from those with $0 \leq w_0 \leq f_0$ of the same order.
4. If $W \in B_F \cup -B_I$ would fulfill $w_0 - w_1 = 2(f_0 - f_1)$, then $V := W - 2F$ would be a basic class orthogonal to $G$, which does not exist. If $W \in -B_I$ with $w_0 - w_1 > 2(f_0 - f_1)$, then also $(W - 2F) \in -B_I$, (in fact of a higher order then $W$). Therefore we obtain all $W \in B_F \cup -B_I$ from those with $0 < w_0 - w_1 < 2(f_0 - f_1)$.

5. If $W \in -B_I$ with $WF \leq 2w_if_i < 0$ for some $i$, then $V := W + 2w_iE_i \in B_F \cup -B_I$. If many of the $f_i$ are 1, this reduces considerably the number of cases to consider; often one can exclude the existence of certain basic classes by excluding that of corresponding elements of $B_F$, which is easy by 3.

To determine the basic classes in the following examples we have used the above observations, and some additional elementary arguments.

(1) If $F = (f_0, f_1, \ldots, f_s, 0, \ldots, 0)$ with $s \leq 8$, then $\Phi_C^{X,F} = 0$ for all $C$ (this follows from Corollary 4.10 and Proposition 4.6).

(2) If $f_1 \geq f_0 - 1$, then $\Phi_C^{X,F} = 0$ for all $C$ (use Remark 5.9, parts 3. and 4.).

(3) If all $f_i$ are odd and $N \geq 9$, then $N - 1$ is divisible by 8 and $(X, F)$ is of strictly $(N-1)/8$-th order simple type. The only element of $B_F \cup -B_I$ of maximal order is $F$. Therefore we get for if $CF$ odd,

$$\Phi_C^{X,F}(e^{xz}(1+p/2)(1-p^2/4)^{(N-9)/8}) = -(-1)^{C(C+F)/2}e^{\mathbf{Q}(x)z^2/2}/\cosh(Fxz)$$

A particular case is that of $N = 9$ and $F = (3, 1^9)$, where $(X, F)$ is of simple type. In this case the result was, for $CF$ odd, brought to our attention by John Morgan, and later Zoltan Szabó informed us of some similar results in other cases (see [M-Sz]). This has been one of the principal motivations for this work.

(4) If $f_0 = 2n$ (with $n \geq 2$) and $\sum_{i \geq 1} f_i \geq 4n^2 - 2n$, and all $f_i$ are strictly positive, then $(X, F)$ is strictly of $\binom{n}{2}$-th order simple type. (We omit the proof.) One basic class of maximal order is $(2n-1, 1^N)$.

(5) If $F = (4, 2^2, 1^8)$ then $(X, F)$ is of simple type (see (4)). $B_F$ consists of $W = (3, 1^{10})$ and $2F - W$ and $B_I = \emptyset$. Therefore we obtain e.g.

$$\Phi_{E_{10}}^{X,F}(e^{xz}(1+p/2)) = e^{\mathbf{Q}(x)z^2/2}\frac{\cosh((H + E_1 + E_2)xz)}{\cosh(Fxz)}$$

(6) If $F = (4, 2, 1^{12})$ then $(X, F)$ is of simple type (see (4)). $B_F \cup -B_I$ consists of

$$W_1 := (3, 1^{13}),\ W_1^i := W_1 - 2E_i\ (i = 2, \ldots, 13),\ W_2^i := (5, 3, 1^{12}) + 2E_i \qquad (i = 2, \ldots, 13).$$

Of these, all but $W_1$ lie in $B_F$. So we obtain, e.g.,

$$\Phi_{E_2}^{X,F}(e^{xz}(1+p/2)) = \frac{e^{Q(x)z^2/2}}{8}\left(2\cosh(W_1xz) + \sum_{i=2}^{13}(1 - 2\delta_{i,2})\frac{\cosh((H + E_1 + 2E_i)xz)}{\cosh(Fxz)}\right).$$

(7) $F = (4, 1^{16})$: Also here $(X, F)$ is of simple type. $B_F \cup -B_I$ consists of

$$W_1 := (3, 1^{16}),\ W_1^i := W_1 - 2E_i\ (1 \leq i \leq 16),$$
$$W_1^{i,j} := W_1 - 2E_i - 2E_j\ (1 \leq i < j \leq 16),\ W_2^{i,j} := W_1 + 2H + 2E_i + 2E_j\ (1 \leq i < j \leq 16).$$

The $W_1^{i,j}$ and the $W_2^{i,j}$ lie in $B_F$.

(8) $F = (5, 2^5, 1^5)$. Also here $(X, F)$ is of simple type. The elements of $B_F \cup -B_I$ are $W = (3, 1^{10})$ and $2F - W$.

(9) If $F = (5, 3, 1^{16})$ then $(X, F)$ is of second order simple type. The unique class in $B_F \cup -B_I$ of order 2 is $F$. The classes of order 1 are

$$W_1 := (3, 1^{17}), \ W_1^i := W_1 - 2E_i \ (2 \leq i \leq 17), \ W_1^{i,j} := W_1 - 2E_i - 2E_j \ (2 \leq i < j \leq 17),$$
$$W_2^i := F + 2E_i \ (2 \leq i \leq 17), W_2^{i,j} := F + 2E_i - 2E_j \ (2 \leq i, j \leq 17, \ i \neq j),$$
$$W_3^{i,j} := F + 2H + 2E_1 + 2E_i + 2E_j \ (2 \leq i < j \leq 17).$$

The classes in $B_F$ are $F$, the $W_1^{i,j}$, the $W_2^{i,j}$ and the $W_3^{i,j}$.

(10) If $F = (5, 1^{25})$ then $(X, F)$ is of 3-rd order simple type. The unique class in $B_F \cup -B_I$ of order 3 is $F$. The classes of order 2 are

$$W_1^i := F + 2E_i, \ (1 \leq i \leq 25), \ W_1^{i,j} := F + 2E_i - 2E_j, (1 \leq i, j \leq 25, \ i \neq j).$$

The classes of order 1 are

$$W_2^J := (3, 1^{25}) - \sum_{j \in J} 2E_j \ (J \subset \{1, \ldots, 25\}, \ |J| \leq 5),$$
$$W_3^{J_1, J_2} := F + \sum_{j_1 \in J_1} 2E_{j_1} - \sum_{j_2 \in J_2} 2E_{j_2} \ (J_1, J_2 \subset \{1, \ldots, 25\} \text{ disjoint}, \ |J_1| = 2, \ |J_2| \leq 2),$$
$$W_4^J := (7, 1^{25}) + \sum_{j \in J} 2E_j, (J \subset \{1, \ldots, 25\}, \ |J| = 5).$$

We finish by giving a list of the types of all boundary points $F$ with $1 \leq f_0 \leq 6$ for which the $\Phi_C^{X,F}$ do not vanish identically. The verifications are elementary, though quite tedious.

1-st order: $(3, 1^9)$, $(4, 1^{16})$, $(4, 2, 1^{12})$, $(4, 2^2, 1^8)$, $(5, 2^3, 1^{13})$, $(5, 2^4, 1^9)$, $(5, 2^5, 1^5)$, $(5, 3, 2, 1^{12})$,

$(5, 3, 2^2, 1^8)$, $(5, 3^2, 1^7)$, $(6, 2^6, 1^{12})$, $(6, 2^7, 1^8)$, $(6, 2^8, 1^4)$, $(6, 3, 2^4, 1^{11})$, $(6, 3, 2^5, 1^7)$,

$(6, 3, 2^6, 1^3)$, $(6, 3^2, 2, 1^{14})$, $(6, 3^2, 2^2, 1^{10})$, $(6, 3^2, 2^3, 1^6)$, $(6, 3^3, 1^9)$.

2-nd order: $(5, 2, 1^{21})$, $(5, 2^2, 1^{17})$, $(5, 3, 1^{16})$, $(6, 2^4, 1^{20})$, $(6, 2^5, 1^{16})$, $(6, 3, 2, 1^{23})$, $(6, 3, 2^2, 1^{19})$,

$(6, 3, 2^3, 1^{15})$, $(6, 3^2, 1^{18})$, $(6, 4, 1^{20})$, $(6, 4, 2, 1^{16})$.

3-rd order: $(5, 1^{24})$, $(6, 1^{36})$, $(6, 2, 1^{32})$, $(6, 2^2, 1^{28})$, $(6, 2^3, 1^{24})$, $(6, 3, 1^{27})$.

## References

[Bo]     Borcherds, R., *Automorphic forms with singularities on Grassmannians*, preprint 1996.

[Br]     Brussee, R., *The canonical class and $C^\infty$ properties of Kähler surfaces*, preprint 1995.

[Do]     Donaldson, S.K., *Polynomial invariants for smooth four-manifolds*, Topology **29** (1990), 257–315.

[E-G1]   Ellingsrud, G. and Göttsche, L., *Variation of moduli spaces and Donaldson invariants under change of polarisation*, J. reine angew. Math. **467** (1995), 1–49.

[E-G2]   Ellingsrud, G. and Göttsche, L., *Wall-crossing formulas, Bott residue formula and the Donaldson invariants of rational surfaces*, preprint 1995.

[E-LP-S] Ellingsrud, G., Le Potier, J. and Strømme, S. A., *Some Donaldson Invariants of $\mathbb{CP}^2$*, preprint 1995.

[E-Z]    Eichler, M., Zagier, D., *The theory of Jacobi forms*, Progress in Math **55** Birkhäuser 1985.

[F-S1]   Fintushel, R., Stern, R.J., *The blowup formula for Donaldson invariants*, Annals of Math. **143** (1996), 529–546.

[F-S2]   Fintushel, R., Stern, R.J., *Donaldson invariants of 4-manifolds of simple type*, J. Diff. Geom. **42** (1995), 577–633.

[F-S3]   Fintushel, R., Stern, R.J., *Rational blowdowns of smooth 4-manifolds*, preprint 1995.

[F-M]    Friedman, R., Morgan, J., *On the diffeomorphism types of certain algebraic surfaces, I*, J. Diff. Geom. **27** (1988), 297–369.

[F-M2]   Friedman, R., Morgan, J., *Algebraic surfaces and Seiberg-Witten invariants*, preprint 1995.

[F-Q]    Friedman, R., Qin, Z., *Flips of moduli spaces and transition formulas for Donaldson polynomial invariants of rational surfaces*, Commun. in Analysis and Geometry **3** (1995), 11–83.

[G]      Göttsche, L., *Modular forms and Donaldson invariants for 4-manifolds with $b_+ = 1$*, J. Amer. Math. Soc. **9** (1996), 827–843.

[H-B-J]  Hirzebruch, F., Berger, T., Jung, R., *Manifolds and Modular forms*, Vieweg, Braunschweig-Wiesbaden, 1992.

[K1]     Kotschick, D., *$SO(3)$-invariants for 4-manifolds with $b^+ = 1$*, Proc. London Math. Soc. **63** (1991), 426–448.

[K2]     Kotschick, D., *On manifolds homeomorphic to $CP^2 \# 8\overline{CP}^2$*, Invent. Math. **95** (1989), 591–600.

[K-L]    Kotschick, D., Lisca, P., *Instanton Invariants of $\mathbb{CP}^2$ via Topology*, Math. Ann. **303** (1995), 345–371.

[K-M]    Kotschick, D., Morgan, J., *$SO(3)$-invariants for 4-manifolds with $b^+ = 1$ II*, J. Diff. Geom. **39** (1994), 433–456.

[Kr-M1]  Kronheimer, P., Mrowka, T., *Recurrence relations and asymptotics for four-manifold invariants*, Bull. Amer. Math. Soc. **30** (1994), 215–221.

[Kr-M2]  Kronheimer, P., Mrowka, T., *Embedded surfaces and the structure of Donaldson's Polynomial invariants*, J. Diff. Geom. **33** (1995), 573-734.

[L-Q]    Li, W.P., Qin, Z., *Lower-degree Donaldson polynomials of rational surfaces*, J. Alg. Geom. **2** (1993), 413–442.

[M-O]    Morgan, J., Ozsváth, P., private communication.

[M-Sz]   Morgan, J., Szabó, Z., Embedded tori in four-manifolds, preprint 1996.

[O-V]    Okoneck, C., Van de Ven, A., *$\Gamma$-type invariants associated to $PU(2)$-bundles and the differentiable structure of Barlow's surface*, Invent. Math. **95** (1989), 601–614.

[O]      Ozsváth, P., *Some Blowup formulas for $SU(2)$ Donaldson Polynomials*, J. Diff. Geom. **40** (1994), 411–447.

[P-T1]   Pidstrigach, V.Y. Tyurin, A.N., *Localization of the Donaldson invariants along Seiberg-Witten classes*, preprint 1995.

[P-T2]   Pidstrigach, V.Y. Tyurin, A.N., in preparation.

[Q1]     Qin, Z., *Equivalence classes of polarizations and moduli spaces of sheaves*, J. Diff. Geom. **37** (1993), 397–413.

[Q2]     Qin, Z., *Moduli of stable sheaves on ruled surfaces and their Picard groups*, J. reine angew. Math. **433** (1992), 201–219.

[R]      Rankin, R.A., *Modular forms and functions*, Cambridge University Press, Cambridge, 1977.

[S-W]    Seiberg, N., Witten, E., *Electric-Magnetic Duality, Monopole Condensation, and Confinement in $N = 2$ Supersymmetric Yang-Mills Theory*, Nucl. Phys. **B279** (1994)

[T1]     Taubes, C.H., in preparation.

[T2]     Taubes, C.H., *The Seiberg-Witten and the Gromov invariants*, Math. Research Letters **2** (1995), 221–238.

[We]     Weil, A., *Elliptic functions according to Eisenstein and Kronecker*, Grundlehren **88**, Springer, Berlin 1976.

[W1]     Witten, E., *Monopoles and Four-Manifolds*, Math. Res. Lett. **1** (1994).

[W2]     Witten, E., *On S-Duality in Abelian Gauge Theory*, Selecta Math. **1** (1995) 383–410.

[Z]      Zagier, D., *Periods of modular forms and Jacobi theta functions*, Invent. math. **104** (1991), 449-465.

INSTITUT MITTAG LEFFLER, AURAVÄGEN 17, S-18262 DJURSHOLM, SWEDEN
*E-mail address*: `gottsche@ml.kva.se`

MAX–PLANCK–INSTITUT FÜR MATHEMATIK, GOTTFRIED–CLAREN–STRASSE 26, D-53225 BONN, GERMANY
*E-mail address*: `zagier@mpim-bonn.mpg.de`